\documentclass[journal]{IEEEtran}

\usepackage[utf8]{inputenc}
\usepackage{graphicx}
\usepackage{cite}       
\graphicspath{ {Figures_bnd_wave/} }
\usepackage{amsmath}
\usepackage{amsfonts}
\usepackage{color}
\usepackage{amssymb}
\usepackage{mathrsfs}
\usepackage{float}
\usepackage{enumerate}
\usepackage{verbatim}
\usepackage{setspace}
\usepackage{epsfig}
\usepackage{url}
\usepackage{graphicx}
\usepackage{lscape}

\newcommand\Tstrut{\rule{0pt}{2.6ex}}
\newcommand\bigZero[1][17]{\mbox{\fontsize{#1}{0}\selectfont$0$}}

\newtheorem{theorem}{Theorem}
\newtheorem{prop}{Proposition}
\newtheorem{corollary}{Corollary}
\newtheorem{lemma}{Lemma}
\newtheorem{remark}{Remark}
\begin{document}

\title{ Finite-dimensional boundary control of the linear Kuramoto-Sivashinsky equation under point measurement with guaranteed $L^2$-gain}		
%

\author{Rami~Katz and
Emilia~Fridman,~\IEEEmembership{Fellow,~IEEE}
\thanks{R. Katz ({\tt\small rami@benis.co.il}) and E. Fridman ({\tt\small  emilia@eng.tau.ac.il}) are with the School of Electrical Engineering, Tel Aviv University, Israel. }
\thanks{Supported by  Israel Science Foundation (grant no. 673/19) and
by Chana and Heinrich Manderman Chair at Tel Aviv University.}%
}
\markboth{}%
{}

\maketitle

\begin{abstract}
Finite-dimensional observer-based controller design for PDEs is a challenging problem.
Recently,  such controllers were introduced for the 1D heat equation, under the assumption that one of the observation or control
operators is bounded.
This paper suggests a constructive method for such controllers
for 1D parabolic PDEs
with both
(observation and control) operators being unbounded. We consider the Kuramoto-Sivashinsky equation (KSE) under either boundary or in-domain point measurement and boundary actuation.
We employ
a modal decomposition approach via dynamic extension, using eigenfunctions of a Sturm-Liouville operator.
The
controller dimension is defined by the number of unstable modes, whereas the observer dimension $N$ may be larger than this number.
We suggest a direct Lyapunov approach to the full-order closed-loop
system, which results in an LMI whose elements and dimension depend on $N$. The value of $N$ and the decay rate are  obtained from the LMI. We  extend our approach to internal stabilization with guaranteed $L^2$-gain  and input-to-state  stabilization 
in the presence of  disturbances in the PDE and the measurement.
We prove that the LMIs are always feasible provided $N$ and the $L^2$ or ISS gains are large enough, thereby obtaining guarantees for our approach. Moreover, for the case of stabilization,
 we show that feasibility of the
 LMI for some $N$ implies its feasibility for $N+1$ (i.e., enlarging $N$ in the LMI cannot deteriorate the resulting decay rate of the closed-loop system).
Numerical examples demonstrate the efficiency of the method.
\end{abstract}
\section{Introduction}
Parabolic PDEs have many applications in physics and engineering. Among such PDEs, the Kuramoto-Sivashinsky equation (KSE) describes many important processes, including chemical reaction-diffusion, flame propagation and viscous flow (see, e.g, \cite{christofides2001,kuramoto1975formation,sivashinsky1977nonlinear,nicolaenko1986some}).

Distributed state-feedback and observer-based control of the KSE was suggested in \cite{armaou2000feedback,christofides2000global} via a modal decomposition approach.
A boundary controller for the KSE in case of a small anti-diffusion parameter was designed in \cite{liu2001stability}. State-feedback stabilization of KSE under boundary or non-local actuation was studied in \cite{cerpa2010null,guzman2019stabilization} by using modal decomposition, whereas null controllability of the KSE was studied in
\cite{cerpa2017control}. 
Stability of the linear KSE as well as its stabilization using a single distributed control were studied in \cite{al2018linearized}.

Output-feedback controllers are more realistic for implementation. Finite-dimensional static output-feedback controllers were suggested in \cite{Aut12,NetzerAut14,lunasin2017finite,
kang2019distributed,selivanov2019delayed} via the spatial decomposition method. However, such controllers may require many sensing and actuation devices.

Observer-based controllers for parabolic equations have been constructed in \cite{curtain1982finite,lasiecka2000control,orlov2004robust,katz2020boundary}, where an observer was designed in the form of a PDE. An advantage of PDE observers is the resulting  separation of controller and observer designs. However, they are often difficult for numerical implementation due to high computational complexity.

Finite-dimensional observer-based controllers for parabolic PDEs were suggested in \cite{curtain1982finite,balas1988finite,christofides2001,harkort2011finite}, whereas finite-dimensional boundary observers for the heat equation were constructed  in \cite{selivanov2018boundary}.
In particular, for bounded control and observation operators, it was shown in \cite{balas1988finite} that the closed-loop system is stable provided the controller dimension is large enough. A singular perturbation approach that reduces the controller design to a finite-dimensional slow system was suggested
in
\cite{christofides2001}, without giving constructive and rigorous conditions for finding the dimension of the slow system that guarantees a desired closed-loop performance of the full-order
system. A bound on the controller dimension  was suggested in \cite{harkort2011finite}. However, this bound was shown to be conservative. Recently
 an efficient bound on the controller dimension in terms of simple LMIs  was  suggested for the 1D heat equation in \cite{RamiContructiveFiniteDim,RamiContructiveFiniteDimDelay} for the case when at least one of the observation or control operators is bounded. The challenging case where both operators are unbounded remained open.

$H_{\infty}$ control of abstract distributed parameter systems was studied in \cite{van2012h}, where the $H_{\infty}$ control problem was reduced to solvability of operator Riccati equations. LMI-based conditions for $H_{\infty}$ control of PDEs and time-delay systems were derived in \cite{NetzerAut14,selivanov2019delayed}, \cite{full} and \cite{Aut09b}. Recently, input-to-state stability (ISS) of PDEs has regained much interest. ISS for the 1D heat equation with boundary disturbance was studied in \cite{Karafyllis2016a}. State-feedback with ISS analysis of diagonal boundary control systems was considered in \cite{lhachemi2020exponential}. Non-coercive Lyapunov functionals for ISS of inifinite-dimensional system were studied in \cite{jacob2019non}. A survey of ISS results can be found in \cite{mironchenko2020input}.

In this paper, for the first time, we provide a constructive method for finite-dimensional observer-based control of a parabolic PDE with the observation and control operators both unbounded. We consider control of the 1D linear KSE under point measurement under either (mixed) Dirichlet or (mixed) Neumann actuation.
This is the first LMI-based method for
finite-dimensional observer-based control of the KSE. We use dynamic extension (see e.g. \cite{curtain}, Sect. 3.3). Dynamic extension was employed for the state-feedback case in \cite{coron2004global,cerpa2010null} and for observer-based control in \cite{RamiContructiveFiniteDimDelay} and
which allows us to manage with unbounded observation and control operators via modal decomposition. Differently from the existing modal decomposition methods for KSE (see, e.g. \cite{cerpa2010null,guzman2019stabilization}), we
introduce a method based on a Sturm-Liouville operator with explicit eigenfunctions and eigenvalues.
 In comparison to \cite{cerpa2010null,guzman2019stabilization},  where the eigenfunctions and eigenvalues can only be approximated numerically,
 our novel approach does not require such approximations.

  We study internal stabilization with guaranteed $L^2$-gain 
  and input-to-state  stabilization in the presence of  disturbances in both the PDE and measurement.
  Note that stabilization with guaranteed $L^2$-gain has not been studied yet via modal decomposition for parabolic PDEs.
In the design, the controller dimension  is defined by the number of unstable modes, whereas the observer's dimension $N$ may be larger than this number. The observer and controller gains are found separately by solving Lyapunov inequalities. We use a direct Lyapunov approach to the full-order closed-loop system.
We derive 
LMIs, whose dimension depends on $N$. These LMIs are used for finding $N$, the resulting exponential decay rate and the $L^2$ and ISS gains.
{We provide feasibility guarantees for the derived LMIs in the cases of $L^2$  and ISS gains
for large enough $N$ and gains.  For the case of stabilization we also prove that feasibility for $N$ implies feasibility for $N+1$ (meaning that the decay rate does not deteriorate).} Numerical examples demonstrate the efficiency of the presented method.

Preliminary results on stabilization of unperturbed 1D KSE under Dirichlet boundary conditions, were presented in \cite{Rami_CDC20}.

\emph{{\bf{Notation:}}}
$ L^2(0,1)$ is the Hilbert space of 
square integrable functions $f:[0,1]\to \mathbb{R} $  with the inner product $\left< f,g\right>:=\scriptsize{\int_0^1 f(x)g(x)dx}$ and induced norm $\left\|f \right\|^2:=\left<f,f \right>$.
$H^{k}(0,1)$ is the Sobolev space of  functions having $k$ square integrable weak derivatives, with the norm $\left\|f \right\|^2_{H^k}:=\sum_{j=0}^{k} \left\|f^{(j)} \right\|^2$. We denote $f\in H^1_0(0,1)$ if $f\in H^1(0,1)$ and $f(0)=f(1)=0$. The Euclidean norm on $\mathbb{R}^n$ is denoted by $\left|\cdot \right|$. For $P \in \mathbb{R}^{n \times n}$, $P>0$ means $P$ is symmetric and positive definite. Sub-diagonal elements of a symmetric matrix are denoted by $*$. For $0<U\in \mathbb{R}^{n\times n}$ and $x\in \mathbb{R}^n$ let $\left|x\right|^2_U:=x^TUx$. $\mathbb{Z}_+$ denotes the nonnegative integers. $\mathbb{N}$ are the natural numbers.
\section{Mathematical preliminaries}
Consider the Sturm-Liouville eigenvalue problem
\begin{equation}\label{eq:SL}
	\begin{aligned}
		\phi''+\lambda \phi = 0,\ \ x\in (0,1)
	\end{aligned}
\end{equation}
with one of the following boundary conditions:
\begin{equation}\label{eq:2BCs}
	\begin{array}{lll}
		&\text{Dirichlet (D): } \phi(0) = \phi(1) = 0,\\
		&\text{Neumann (Ne): } \phi'(0) = \phi'(1)=0.
	\end{array}
\end{equation}
These problems induce a sequence of eigenvalues $\lambda_n$ with corresponding eigenfunctions $\phi_n^D(x)$ and $\phi_n^{Ne}(x)$ given by
\begin{equation}\label{eq:SLBCs}
\begin{array}{lll}
&\text{(D): } \lambda_n = n^2\pi^2,\quad  \phi_n^D(x)=\sqrt{2}\sin\left(\sqrt{\lambda_n} x\right),\ n\in \mathbb{N},\\
&\text{(Ne): } \lambda_0=0,\quad \quad  \lambda_n = n^2\pi^2,\\
&\hspace{0.77cm} \phi_0^{Ne}(x)\equiv 1, \quad \phi_n^{Ne}(x)=\sqrt{2}\cos\left(\sqrt{\lambda_n}x\right),\ n\in \mathbb{N}.
\end{array}
\end{equation}
The eigenfunctions form complete and orthonormal family in $L^2(0,1)$.
\begin{lemma}\label{lem:H1Equiv}\cite{renardy2006introduction}
Let $h \overset{L^2}{=} \sum_{n=1}^{\infty}h_n\phi_n^D$. 	Then $h\in H^1_0(0,1)$ if and only if $\sum_{n=1}^{\infty}\lambda_nh_n^2< \infty$. Moreover, \begin{equation}\label{lem22}
\left\|h' \right\|^2=\sum_{n=1}^{\infty}\lambda_nh_n^2.
\end{equation}
\end{lemma}
\begin{lemma}\label{lem:H1Equiv1}
Assume $h \overset{L^2}{=} \sum_{n=0}^{\infty}h_n\phi_n^{Ne}$. Then $h\in H^2(0,1)$ with $h'(0)=h'(1)=0$ if and only if $\sum_{n=1}^{\infty}\lambda_n^2h_n^2< \infty$.  Moreover,
\begin{equation}\label{eq:H2}
\left\|h'' \right\|^2=\sum_{n=1}^{\infty}\lambda_n^2h_n^2.
\end{equation}
\end{lemma}
\begin{IEEEproof}
Assume $h\in H^2(0,1)$ with $h'(0)=h'(1)=0$. Let $n\in \mathbb{N}$. Integrating by parts twice and taking into account \eqref{eq:2BCs} we have $-\left<h'',\phi_n^{Ne} \right>= \lambda_n\left<h,\phi_n^{{Ne}} \right>$. Applying Parseval's equality, we have $\sum_{n=1}^{\infty}\lambda_n^2h_n^2< \infty$. For the other direction, assume $\sum_{n=1}^{\infty}\lambda_n^2h_n^2< \infty$. Given $N\in \mathbb{N}$, let
\begin{equation*}
T_N(x) = \sum_{n=0 }^Nh_n\phi_n^{N_e}(x), \ S_N(x) = -\sum_{n=1 }^N\lambda_nh_n\phi_n^{N_e}(x).
\end{equation*}
By assumption, $\left\{S_N\right\}_{N\in \mathbb{N}}$ converge in $L^2(0,1)$ to $S = -\sum_{n=1}^{\infty}\lambda_nh_n\phi_n^{N_e}$. Take any smooth function $\rho(x)$, compactly supported in $(0,1)$. Then, integration by parts gives $\left<T_N,\rho'' \right> = \left<S_N,\rho \right>$. Since $T_N \overset{N\rightarrow \infty}{\longrightarrow} h$ in $L^2(0,1)$, taking $N\to \infty$ we obtain $\left<h,\rho'' \right> = \left<S,\rho \right>$. Thus, $h$ has a weak derivative of second order $h'' = S \in L^2(0,1)$. We deduce $h\in H^2(0,1)$. In particular, by Sobolev's embedding theorem, $h\in C^1(0,1)$. Furthremore, boundedness of $\left\{\phi_n^{Ne}\right\}_{ n\in \mathbb{Z}_+}$ on $[0,1]$ and
\begin{equation*}\label{eq:term by term}
\sum_{n=1}^{\infty}\left|h_n \right|\leq \left(\sum_{n=1}^{\infty}\lambda_n^2h_n^2 \right)^{\frac{1}{2}} \left(\sum_{n=1}^{\infty}\lambda_n^{-2} \right)^{\frac{1}{2}}<\infty
\end{equation*}
imply uniform convergence of the series and $h \overset{C}{=} \sum_{n=0}^{\infty}h_n\phi_n^{Ne}$ (continuous functions which agree almost everywhere). Since $\left(\phi_n^{Ne}\right)'(x) = -\sqrt{\lambda_n}\phi_n^D(x)$ for $n\in \mathbb{N}$, with $\left\{\phi_n^D\right\}_{n\in \mathbb{N}}$ bounded on $[0,1]$,
\begin{equation*}\label{eq:term by term1}
\sum_{n=1}^{\infty}\sqrt{\lambda_n}\left|h_n \right|\leq \left(\sum_{n=1}^{\infty}\lambda_n^2h_n^2 \right)^{\frac{1}{2}} \left(\sum_{n=1}^{\infty}\lambda_n^{-1} \right)^{\frac{1}{2}}<\infty
\end{equation*}
shows that the series can be differentiated term-by-term. Differentiating term by term we get
\begin{equation*}
h' = -\sum_{n=1}^{\infty}\sqrt{\lambda_n}h_n\phi_n^{D}\Longrightarrow \left\|h'\right\|^2 = \sum_{n=1}^{\infty}\lambda_nh_n^2,
\end{equation*}
where the right-hand side follows by orthonormality of $\left\{\phi_n^D\right\}_{n=1}^{\infty}$. Finally, substituting $x\in \left\{0,1\right\}$ into $h'$, we have $h'(0)=h'(1)$.
\end{IEEEproof}
\begin{lemma}\label{lem_sob}
(Sobolev's inequality \cite{kang2019distributed})
Let	$h\in H^1(0,1)$. Then, for all $\Gamma>0$
\begin{equation*}
\max_{x\in[0,1]}\left|h(x) \right|^2\leq (1+\Gamma)\left\|h\right\|^2+ \Gamma^{-1}\left\|h' \right\|^2.
\end{equation*}
\end{lemma}
\begin{remark}
Differently from Theorem 8.8 in \cite{brezis2010functional}, Lemma \ref{lem_sob} gives an explicit upper bound on $\left\|h\right\|_{L^{\infty}(0,1)}$, depending on a general constant $\Gamma>0$. A variant of Lemma \ref{lem_sob} with $\Gamma = 1$ was given in \cite{zheng2018input}.
\end{remark}

\section{Stabilization of the linear 1D KSE}
In this section we consider stabilization of the linear 1D Kuramoto-Sivashinsky equation (KSE)
\begin{equation}\label{eq:LinearizedKSE}
	z_t(x,t) = -z_{xxxx}(x,t) -\nu z_{xx}(x,t),
\end{equation}
where $t\geq 0$, $x\in (0,1)$, $z(x,t)\in \mathbb{R}$ and $\nu>0$ is the "anti-diffusion" coefficient.

We consider either (mixed) Dirichlet boundary conditions
\begin{equation}\label{eq:BCSD}
	\begin{array}{lll}
		&\text{(D)}\ z(0,t) = u(t),  \quad z(1,t)= 0 ,\\
		&\hspace{6mm}z_{xx}(0,t)=0, \quad z_{xx}(1,t)=0
	\end{array}
\end{equation}
or (mixed) Neumann boundary conditions
\begin{equation}\label{eq:BCSNe}
	\begin{array}{lll}
		&\text{(Ne)}\ z_x(0,t) = u(t),  \quad z_x(1,t)= 0 ,\\
		&\hspace{8mm}z_{xxx}(0,t)=0, \quad z_{xxx}(1,t)=0.
	\end{array}
\end{equation}
For both cases $u(t)$ is a control input to be designed.

The boundary conditions \eqref{eq:BCSD} and \eqref{eq:BCSNe} have been considered in \cite{anders2012higher}. A detailed description of KSE with either \eqref{eq:BCSD} or \eqref{eq:BCSNe} can be found in \cite{nicolaenko1986some}. The boundary conditions \eqref{eq:BCSD} and \eqref{eq:BCSNe} allow to use modal decomposition with respect to the eigenfunctions of \eqref{eq:SL} and \eqref{eq:2BCs} in order to obtain either $H^1(0,1)$ (Dirichlet) or $H^2(0,1)$ (Neumann) stability of the closed-loop system. See also Remark \ref{rem_bc} below about modal decomposition approach under other boundary conditions.\\[0.2cm]
\emph{A. Dirichlet actuation and in-domain point measurement}\label{Sec:1}

Consider the KSE \eqref{eq:LinearizedKSE} with boundary conditions \eqref{eq:BCSD} and in-domain point measurement
\begin{equation}\label{eq:InDomPointMeas}
	y(t) = z(x_*,t), \ x_*\in (0,1).
\end{equation}
We introduce the change of variables
\begin{equation}\label{eq:ChangeVars}
	w(x,t)=z(x,t)-r(x)u(t), \quad r(x):=1-x
\end{equation}
to obtain the following equivalent
ODE-PDE system
\begin{equation}\label{eq:PDE1PointActConstChangeVars}
	\begin{array}{lll}
		& \dot{u}(t)=v(t),\\
		& w_t(x,t)=-w_{xxxx}(x,t)-\nu w_{xx}(x,t)-r(x)v(t)
	\end{array}
\end{equation}
with 
boundary conditions
\begin{equation}\label{eq:BCsDirChangeVars}
\begin{array}{lll}
& w(0,t)=0, \quad w(1,t)=0,\\
& w_{xx}(0,t) = 0, \quad w_{xx}(1,t)=0.
\end{array}
\end{equation}
and measurement
\begin{equation}\label{eq:InDomPointMeas1}
	y(t) = w(x_*,t)+r(x_*)u(t).
\end{equation}
Henceforth we treat $u(t)$ as an additional state variable and $v(t)$ as the control input. Given $v(t)$, $u(t)$ can be computed by integrating $\dot{u}(t)=v(t)$, where we choose  $u(0)=0$.

Differently from state-feedback control (see e.g. \cite{mironchenko2020local}), our  output-feedback control law  will be coupled with the PDE through the measurement \eqref{eq:InDomPointMeas1}. Therefore, for well-posedness, we will consider the closed-loop system consisting of \eqref{eq:PDE1PointActConstChangeVars} and the ODEs \eqref{eq:WobsODENonDelayed}, which define the control input \eqref{eq:WContDef} (see \eqref{eq:DiffOpKSE}-\eqref{eq:ClassicalSolDirichlet2} below). We will show that the closed-loop system, subject to the proposed control law \eqref{eq:WContDef}, has a unique classical solution. Thus, the use of modal decomposition in \eqref{eq:Wseries} and \eqref{eq:WOdesNonDelayed} below will be justified a posteriori and is presented here in order to  construct a finite-dimensional observer-based controller.

We present the solution to \eqref{eq:PDE1PointActConstChangeVars} as
\begin{equation}\label{eq:Wseries}
	\begin{array}{lll}
		w(x,t) &= \sum_{n=1}^{\infty}w_n(t)\phi_n^D(x), \ 	w_n(t) =\left<w(\cdot,t),\phi_n^D\right>
	\end{array}
\end{equation}
with $\left\{\phi_n^D\right\}_{n\in \mathbb{N}}$ defined in \eqref{eq:SLBCs}. Differentiating under the integral sign, integrating by parts and using \eqref{eq:SL} and \eqref{eq:2BCs} we have
\begin{equation}\label{eq:WOdesNonDelayed}
	\begin{array}{lll}
		&\dot{w}_n(t) = (-\lambda_n^2+\nu\lambda_n)w_n(t) + b_nv(t), \\
		&w_n(0)=\left<w(\cdot,0),\phi_n^D\right>, \ b_n=-\left<r,\phi_n^D\right>=-\sqrt{\frac{2}{\lambda_n}}.
	\end{array}
\end{equation}
In particular, note that
\begin{equation}\label{eq:AssbnNonDelayed}
	\begin{array}{lll}
		&b_n \neq 0 , \quad n\geq 1
	\end{array}
\end{equation}
and
\begin{equation}\label{eq:AssbnNonDelayed1}
	\sum_{n=N+1}^{\infty}b_n^2 \leq \frac{2}{\pi^2}\int_N^{\infty}\frac{dx}{x^2}=\frac{2}{\pi^2N}, \ N\geq 1.
\end{equation}
Let $\delta>0$ be a desired decay rate. Since $\lim_{n \to \infty}\lambda_n=\infty$, there exists some $N_0 \in \mathbb{N}$ such that
\begin{equation}\label{eq:N0}
	-\lambda_n^2+\nu \lambda_n<-\delta, \quad n>N_0.
\end{equation}
Let $N\in \mathbb{N}, \ N_0\leq N$. $N_0$ will define the dimension of the controller, whereas $N$ will be the dimension of the observer.

We construct a finite-dimensional observer of the form
\begin{equation}\label{eq:WhatSeries}
	\hat{w}(x,t) = \sum_{n=1}^{N}\hat{w}_n(t)\phi_n^D(x),
\end{equation}
where $\hat{w}_n(t)$ satisfy the ODEs
\begin{equation}\label{eq:WobsODENonDelayed}
	\begin{array}{lll}
		\dot{\hat{w}}_n(t) &= (-\lambda_n^2+\nu \lambda_n)\hat{w}_n(t) + b_nv(t)\\
		&-l_n\left[\hat{w}(x_*,t)+r(x_*)u(t) - y(t)\right],\\
		\hat{w}_n(0)&=0, \quad 1\leq n\leq N.
	\end{array}
\end{equation}
with $y(t)$ given in \eqref{eq:InDomPointMeas1} and scalar observer gains $\left\{l_n\right\}_{n=1}^N$.

We assume the following:

\textbf{Assumption 1:} The point $x_*\in (0,1)$ satisfies
\begin{equation}\label{eq:AsscnNonDelayed}
	c_n=\phi_n^D(x_*)=\sqrt{2}\sin\left(\sqrt{\lambda_n} x_*\right)\neq 0 , \ 1\leq n \leq N_0.
\end{equation}
Assumption 1 is satisfied for $N_0=1$ by any $x^*\in(0,1)$, whereas for $N_0>1$ the corresponding $x^*$ is subject to the following condition:
$x^*\neq k/n<1,\ k=1,..., N_0-1,\ n=2,...,N_0$. E.g, for $N_0=2$ the condition is $x^*\neq 0.5$.

\textbf{Assumption 2:} Assume $$\nu \notin \left\{\pi^2(n^2+m^2) \ ;\ n,m\geq 0 \ n\neq m\right\}\cup\left\{0\right\}.$$
Denote
\begin{equation}\label{eq:C0A0}
\begin{array}{lll}
&A_0 = \operatorname{diag}\left\{-\lambda_1^2+\nu \lambda_1,\dots,-\lambda_{N_0}^2+\nu \lambda_{N_0} \right\},\\
& C_0=\left[c_1,\dots,c_{N_0} \right],\ \tilde{B}_0= \left[1,b_1,\dots,b_{N_0} \right]^T,\\
&\tilde{A}_0 =\operatorname{diag}\left\{0,A_0\right\}\in \mathbb{R}^{(N_0+1)\times(N_0+1)}.
\end{array}
\end{equation}
Under Assumptions 1 and 2 the pair $(A_0,C_0)$ is observable, by the Hautus lemma. We choose $L_0 = \left[l_1,\dots,l_{N_0} \right]^T$ which
 satisfies the Lyapunov inequality
\begin{equation}\label{eq:GainsDesignL}
	P_{\text{o}}(A_0-L_0C_0)+(A_0-L_0C_0)^TP_{\text{o}} < -2\delta P_{\text{o}}
\end{equation}
with $0<P_{\text{o}}\in \mathbb{R}^{N_0\times N_0}$.  
Furthermore, let $l_n=0$ for $n>N_0$.\\[0.1cm]
Assumption 2 and \eqref{eq:AssbnNonDelayed} imply that the pair $(\tilde{A}_0,\tilde{B}_0)$ is controllable, by the Hautus lemma (see also Lemma 6 in \cite{guzman2019stabilization}, where the Kalman rank condition is used). Let $K_0\in \mathbb{R}^{1\times (N_0+1)}$ satisfy
\begin{equation}\label{eq:GainsDesignK}
	\begin{aligned}
		&P_{\text{c}}(\tilde{A}_0+\tilde{B}_0K_0)+(\tilde{A}_0+\tilde{B}_0K_0)^TP_{\text{c}} < -2\delta P_{\text{c}},
	\end{aligned}
\end{equation}
with $0<P_{\text{c}}\in \mathbb{R}^{(N_0+1)\times (N_0+1)}$.

We propose a $(N_0+1)$-dimensional controller of the form
\begin{equation}\label{eq:WContDef}
\begin{aligned}
&\hspace{-2mm}v(t)= K_0\hat{w}^{N_0}(t),\  \hat{w}^{N_0}(t) = \left[u(t),\hat{w}_1(t),\dots,\hat{w}_{N_0}(t) \right]^T
\end{aligned}
\end{equation}
which is based on the $N$-dimensional observer \eqref{eq:WobsODENonDelayed}.


For well-posedness of the closed-loop system \eqref{eq:PDE1PointActConstChangeVars}, \eqref{eq:WobsODENonDelayed} and \eqref{eq:WContDef} we consider the operator
\begin{equation}\label{eq:DiffOpKSE}
	\begin{array}{lll}
		&\mathcal{A}:\mathcal{D}(\mathcal{A})\to L^2(0,1),\ \mathcal{A} = \partial_{xxxx}+\nu \partial_{xx},
	\end{array}
\end{equation}
where
\begin{equation}\label{eq:D1}
	\mathcal{D}(\mathcal{A}) = \left\{h\in H^4(0,1) | h(0)=h(1)=h''(0)=h''(1)=0 \right\}
\end{equation}
is dense in $L^2(0,1)$. Let $h\in \mathcal{D}(\mathcal{A})$. It can be shown using integration by parts that
\begin{equation}\label{eq:CalASeries}
	\mathcal{A}h \overset{L^2}{=} \sum_{n=1}^{\infty}\left(\lambda_n^2-\nu \lambda_n\right)\left<h,\phi_n^D \right>\phi_n^D.
\end{equation}
Furthermore, $\left\{\phi_n^D\right\}_{n\in \mathbb{N}}$ is a complete family of orthonormal eigenfunctions of $\mathcal{A}$. Thus, by Section 2.6 in \cite{GeorgeBook}, $-\mathcal{A}$ is a diagonalizable operator. By Remark 2.6.4 in \cite{GeorgeBook}, $\text{spec}(-\mathcal{A})=\left\{-\lambda_n^2+\nu \lambda_n\right\}_{n=1}^{\infty}$. Since $\lambda_n \to \infty$ as $n\to \infty$, the resolvent set $\rho \left(-\mathcal{A} \right)$ contains a half plane $\left\{z\in \mathbb{C} \ | \ \operatorname{Re}(z)>w\right\}$ for large enough $w\in \mathbb{R}$. Therefore, $-\mathcal{A}$ is a sectorial operator which generates an analytic semigroup on $L^2(0,1)$ (see also Theorem 12.31 in \cite{renardy2006introduction}).

Let $\mathcal{H}:=L^2(0,1)\times \mathbb{R}^{N+1}$ be a Hilbert space with the norm $\left\|\cdot\right\|_{\mathcal{H}}^2:=\left\| \cdot\right\|^2+\left|\cdot \right|^2$. Introducing the state
\begin{equation*}
	\begin{array}{lll}
		&\xi(t) = \text{col}\left\{\xi^{(1)},\xi^{(2)}\right\},\\
		&\xi^{(1)} = w(\cdot,t),\ \ \xi^{(2)} = \text{col}\left\{u(t),\hat{w}_1(t),\dots, \hat{w}_N(t) \right\},
	\end{array}
\end{equation*}
the closed-loop system can be presented as
\begin{equation}\label{eq:AbstractODEDirichlet}
	\begin{array}{lll}
		&\frac{d\xi (t)}{dt}+\tilde{\mathcal{A}}\xi(t)=F(\xi),
\end{array}
\end{equation}
where
\begin{equation}\label{A1}
\begin{array}{lll}
&\tilde{\mathcal{A}}:=diag \{\mathcal{A},\mathcal{B}\}, \ \mathcal{B} = \begin{bmatrix}-\tilde{A}_0-\tilde{B}_0K_0+L_0C_0 & L_0C_1\\
-B_1K_0 & -A_1 \end{bmatrix},\\
& A_1 = \operatorname{diag}\left\{-\lambda_{N_0+1}^2+\nu \lambda_{N_0+1},\dots,-\lambda_{N}^2+\nu \lambda_{N} \right\},\\
&B_1 =\left[b_{N_0+1},\dots,b_N \right]^T,\ C_1 =\left[c_{N_0+1},\dots, c_N \right],\\
&F(\xi) = \begin{bmatrix}f_1(\xi)\\f_2(\xi) \end{bmatrix},\ f_1(\xi)=-r(x)K_0\hat{w}^{N_0}(t),\\
&f_2(\xi)=\text{col}\left\{L_0w(x_*,t),0\right\}.\\
\end{array}
\end{equation}
Here $-\tilde{\mathcal{A}}$ generates an analytic semigroup (since $-\mathcal{A}$ generates an analytic semigroup on $L^2(0,1)$ and $\mathcal{B}$ is a linear operator on $\mathbb{R}^{N+1}$) on $\mathcal{H}$ and the function $F:\mathcal{D}(\mathcal{A})\times \mathbb{R}^{N+1}\to \mathcal{H}$ is linear. Moreover, since for any $h\in \mathcal{D}(\mathcal{A})$
\begin{equation*}
	h(x_*)=\int_0^{x_*}h'(x)dx,
\end{equation*}
we obtain for $\xi\in \mathcal{D}(\mathcal{A})\times \mathbb{R}^{N+1}$
\begin{equation*}
\begin{array}{lll}
&\left\|f_1(\xi)\right\|^2\leq \left\|r\right\|^2\cdot \left|K_0\right|^2\cdot \left\| \xi\right\|_{\mathcal{H}}^2,\\[0.5mm]
&\left|f_2(\xi) \right|^2\leq \text{const}_1\cdot \left\|w_x(\cdot,t) \right\|^2\\[0.5mm]
&\overset{\eqref{lem22}}{=} \text{const}_1 \cdot \sum_{n=1}^{\infty}\lambda_n\left|\left<w(\cdot,t),\phi_n^D\right> \right|^2 \overset{\eqref{eq:CalASeries}}{\leq} \cdot\left[\left\|\xi\right\|_{\mathcal{H}}^2+\left\|\tilde{\mathcal{A}}\xi\right\|_{\mathcal{H}}^2\right].
\end{array}
\end{equation*}
By Theorems 6.3.1 and 6.3.3 in \cite{pazy1983semigroups}, the system \eqref{eq:PDE1PointActConstChangeVars}, \eqref{eq:WobsODENonDelayed} with control input \eqref{eq:WContDef} and initial condition $w(\cdot,0)\in \mathcal{D}(\mathcal{A})$ has a unique classical solution
\begin{equation}\label{eq:ClassicalSolDirichlet1}
\begin{array}{lll}
&\xi \in C([0,\infty);\mathcal{D}\left(\mathcal{A}\right))\cap C^1([0,\infty);\mathcal{H})
\end{array}
\end{equation}
such that
\begin{equation}\label{eq:ClassicalSolDirichlet2}
\xi(t)\in \mathcal{D}(\mathcal{A})\times \mathbb{R}^{N+1}, \ \  t>0.
\end{equation}
Let
\begin{equation}\label{eq:WEstErrorNonDelayed}
e_n(t) = w_n(t)-\hat{w}_n(t), \ 1\leq n \leq N
\end{equation}
be the estimation error. By using \eqref{eq:Wseries} and \eqref{eq:WhatSeries}, the innovation term $\hat{w}(x_*,t)+r(x_*)u(t) - y(t)$ in \eqref{eq:WobsODENonDelayed} can be presented as
\begin{equation}\label{eq:WIntroZetaNonDelayed}
\begin{array}{ll}
&\hspace{-2mm}\hat{w}(x_*,t)+r(x_*)u(t) - y(t) =-\sum_{n=1}^{N} c_ne_n(t)-\zeta_N(t)
\end{array}
\end{equation}
where
\begin{equation}\label{eq:zetaintegral}
\begin{array}{ll}
&\zeta_N(t) = w(x_*,t)-\sum_{n=1}^N w_n(t)\phi_n^D(x_*)\\
&=\int_0^{x_*}\left[w_x(x,t)-\sum_{n=1}^Nw_n(t)\frac{d}{dx}\phi_n^D(x) \right]dx.
\end{array}
\end{equation}
Then the error equations have the form
\begin{equation}\label{eq:Wen}
\begin{array}{ll}
&\dot e_n(t)=(-\lambda_n^2+\nu \lambda_n)e_n(t)\\
&\hspace{6mm}-l_n\left(\sum_{n=1}^{N} c_ne_n(t)+\zeta_N(t)\right),  \  1\leq n \leq N_0,\\
& \dot{e}_n(t)=(-\lambda_n^2+\nu \lambda_n)e_n(t), \ \ N_0+1\leq n \leq N.
\end{array}
\end{equation}
Note that the Cauchy-Schwarz inequality implies
\begin{equation}\label{eq:WzetaEst}
\begin{array}{lll}
&\hspace{-3mm}\zeta_N^2(t)\leq \left(\int_0^{x_*}\left|w_x(x,t)-\sum_{n=1}^Nw_n(t)\frac{d}{dx}\phi_n^D(x) \right|dx \right)^2\\
&\hspace{-3mm}\leq \left\|w_x(\cdot,t)-\sum_{n=1}^Nw_n(t)\frac{d}{dx}\phi_n^D(\cdot) \right\|^2\overset{\eqref{lem22}}{=} \sum_{n=N+1}^{\infty}\lambda_nw_n^2(t).
\end{array}
\end{equation}
Denoting
\begin{equation}\label{eq:ErrDefNonDelayed}
\begin{array}{lllllll}
& X_{N}(t) = \text{col}\left\{\hat{w}^{N_0}(t),e^{N_0}(t),\hat{w}^{N-N_0}(t), e^{N-N_0}(t) \right\},\\
&e^{N_0}(t)=\left[e_1(t),\dots,e_{N_0}(t) \right],\\
&e^{N-N_0}(t)=\left[e_{N_0+1}(t),\dots,e_{N}(t) \right]^T,\\
&\hat{w}^{N-N_0}(t)=\left[\hat{w}_{N_0+1}(t),\dots,\hat{w}_{N}(t) \right]^T,
\end{array}
\end{equation}
and using
\eqref{eq:WOdesNonDelayed}, \eqref{eq:WobsODENonDelayed}, \eqref{eq:WContDef} we arrive at the closed-loop system
\begin{equation}\label{eq:ClosedLoopNonDelayed}
	\begin{aligned}
		&\dot{X}_N(t) = FX_N(t)+\mathcal{L}\zeta_N(t),\quad t\ge 0,\\
		& \dot{w}_n(t) = (-\lambda_n^2+\nu\lambda_n)w_n(t) +b_n\tilde{K}_0X_N(t), \ \ n>N.
	\end{aligned}
\end{equation}
Here
\begin{equation}\label{eq:ErrDefNonDelayed1}
\begin{array}{lllllll}
& \hspace{-1.5mm} F = \scriptsize\begin{bmatrix}\tilde{A}_0+\tilde{B}_0K_0 & \tilde{L}_0C_0 & 0 &\tilde{L}_0C_1 \\ 0 & A_0-L_0C_0 & 0 & -L_0C_1\\ B_1K_0 & 0 & A_1 & 0\\
0 & 0 & 0 & A_1 \end{bmatrix}, \tilde{K}_0 = \begin{bmatrix} K_0,&0_{1\times (2N-N_0)}\end{bmatrix},\\
&\hspace{-1.5mm} \tilde{L}_0 = \text{col}\left\{0,L_0 \right\}\in \mathbb{R}^{N_0+1}, \ \mathcal{L}= \text{col}\left\{\tilde{L}_0,-L_0,0\right\}\in \mathbb{R}^{2N+1}.
	\end{array}
\end{equation}
For stability analysis of the closed-loop system \eqref{eq:ClosedLoopNonDelayed} we consider the Lyapunov function
\begin{equation}\label{eq:VNonDelayed}
V(t)=\left|X_N(t)\right|^2_P+\sum_{n=N+1}^{\infty}\lambda_n w_n^2(t),
\end{equation}
where $0<P\in \mathbb{R}^{(2N+1)\times (2N+1)}$. This Lyapunov function is chosen to compensate $\zeta_N(t)$ using \eqref{eq:WzetaEst}. To justify differentiation of the series in \eqref{eq:VNonDelayed} term-by-term it is sufficient to show that the series of term-by-term derivatives converges uniformly on compact subsets of $(0,\infty)$. Since $\lambda_n\to \infty$ as $n\to \infty$, this reduces to showing that $\sum_{n=1}^{\infty}\lambda_n^3w_n^2(t)$ converges uniformly on compact subsets of $(0,\infty)$. Recall that we have a classical solution satisfying \eqref{eq:ClassicalSolDirichlet1} and \eqref{eq:ClassicalSolDirichlet2}. From \eqref{eq:CalASeries} we find that $\left\|\mathcal{A}w(\cdot,t)\right\|^2 = \sum_{n=1}^{\infty}\left( \lambda_n^2-\nu\lambda_n\right)^2w_n^2(t)$  is continuous on $(0,\infty)$.
Since $\lambda_n\to \infty$ as $n\to \infty$, we have $\sum_{n=1}^{\infty}\lambda_n^3w_n^2(t)\leq M \left\|\mathcal{A}w(\cdot,t)\right\|^2 $ for some constant $M>0$, independent of $t$. Thus $\sum_{n=1}^{\infty}\lambda_n^3w_n^2(t)$ is uniformly bounded on compact sets in $(0,\infty)$. To show uniform convergence, we apply Dini's theorem. Indeed, $\Sigma_N(t) = \sum_{n=1}^{N} \lambda_n^3 w_n^2(t)$ is a sequence of monotonically increasing continuous functions converging pointwise to $\Sigma(t) = \sum_{n=1}^{\infty} \lambda_n^3 w_n^2(t)$. Let $t_1,t_2\in J \subseteq (0,\infty)$, where $J$ is compact. Then,
\begin{equation*}
\begin{array}{lll}
&\left|\sum_{n=1}^{\infty} \lambda_n^3 w_n^2(t_1)- \sum_{n=1}^{\infty} \lambda_n^3 w_n^2(t_2)\right|\\
&\leq \sum_{n=1}^{\infty}\left|\lambda_n^{1.5}\left(w_n(t_1)-w_n(t_2) \right) \right|\left|\lambda_n^{1.5}\left(w_n(t_1)+w_n(t_2) \right) \right|\\
&\leq \left[\sum_{n=1}^{\infty}\lambda_n^3\left(w_n(t_1)+w_n(t_2) \right)^2 \right]^{0.5}\\
&\hspace{15mm}\times\left[\sum_{n=1}^{\infty}\lambda_n^3\left(w_n(t_1)-w_n(t_2) \right)^2 \right]^{0.5}.
\end{array}
\end{equation*}
We now get
\begin{equation*}
\begin{array}{lll}
& \sum_{n=1}^{\infty}\lambda_n^3\left(w_n(t_1)+w_n(t_2) \right)^2 \\
&\hspace{8mm}\leq 2\sum_{n=1}^{\infty}\lambda_n^3\left[w_n^2(t_1)+w_n^2(t_2) \right]\leq \operatorname{const} , \ t_1,t_2\in J,\\
&\sum_{n=1}^{\infty}\lambda_n^3\left(w_n(t_1)-w_n(t_2) \right)^2\\
&\hspace{8mm}\leq M \Vert \mathcal{A} \left(w(\cdot,t_1)-w(\cdot,t_2) \right)\Vert^2\to 0, \ \text{as} \ t_2\to t_1,
\end{array}
\end{equation*}
where the upper was shown and the lower follows since we work with a classical solution. Hence $\Sigma(t) = \sum_{n=1}^{\infty} \lambda_n^3 w_n^2(t)$ is continuous and $\Sigma_N(t)$ converge to $\Sigma(t)$ uniformly. Differentiation of $V(t)$ along the solution of \eqref{eq:ClosedLoopNonDelayed} gives
\begin{equation}\label{eq:WStabAnalysisNonDelayed}
\begin{array}{lll}
&\dot{V}+2\delta V = X_N^T(t)\left[PF +F^TP+2\delta P\right]X_N(t)\\[0.5mm]
&+2X_N^T(t)P\mathcal{L}\zeta_N(t) +2\sum_{n=N+1}^{\infty}(-\lambda_n^3+\nu \lambda_n^2+\delta\lambda_n)w_n^2(t)\\[0.5mm] &+2\sum_{n=N+1}^{\infty}\lambda_n w_n(t)b_n\tilde{K}_0X_N(t).
\end{array}
\end{equation}
The Cauchy-Schwarz inequality implies
\begin{equation}\label{eq:WCrosTermNonDelayed}
	\begin{array}{lllll}
		&\sum_{n=N+1}^{\infty} 2\lambda_n w_n(t)b_n\tilde{K}_0X_N(t)\\
		&\leq \frac{1}{\alpha} \sum_{n=N+1}^{\infty}\lambda_n^2 w_n^2(t)+\alpha\left[\sum_{n=N+1}^{\infty}b_n^2 \right]\left|\tilde{K}_0X_N(t)\right|^2\\
		&\overset{\eqref{eq:AssbnNonDelayed1}}{\leq} \frac{1}{\alpha} \sum_{n=N+1}^{\infty}\lambda_n^2 w_n^2(t)+ \frac{2\alpha }{\pi^2N}\left|\tilde{K}_0X_N(t)\right|^2
	\end{array}
\end{equation}
where $\alpha>0$. From  monotonicity of $\lambda_n, \ n\in \mathbb{N}$ we have
\begin{equation}\label{eq:TailEstDir}
\begin{array}{lll}
&\hspace{-5mm}\sum_{n=N+1}^{\infty}(-\lambda_n^3+\nu \lambda_n^2+\delta\lambda_n+\frac{\lambda_n^2}{2\alpha})w_n^2(t)\\
&\leq -2\left(\theta_{N+1}^{(1)}-
\frac{\lambda_{N+1}}{2\alpha} \right)\sum_{n=N+1}^{\infty}\lambda_nw_n^2(t)\\
&\overset{\eqref{eq:WzetaEst}}{\leq} -2\left(\theta_{N+1}^{(1)}-\frac{\lambda_{N+1}}{2\alpha}\right)\zeta_N^2(t),\\
&\hspace{-5mm}\theta_n^{(1)}=\lambda_n^2-\nu \lambda_n-\delta,\ \ n\geq 1
\end{array}
\end{equation}
if $-\theta_{N+1}^{(1)}+\frac{\lambda_{N+1}}{2\alpha}\leq 0$. Let $\eta(t) = \text{col}\left\{X_N(t),\zeta_N(t) \right\}$. From \eqref{eq:WStabAnalysisNonDelayed}, \eqref{eq:WCrosTermNonDelayed} and  \eqref{eq:TailEstDir} we obtain
\begin{equation}\label{eq:WStabResultNonDelayed}
\begin{array}{ll}
&\dot{V}+2\delta V \leq \eta^T(t)\Psi^{(1)}_N \eta(t)\leq0
\end{array}
\end{equation}
provided
\begin{equation}\label{eq:WLMIsNonDelayed}
\begin{aligned}
&\Psi^{(1)}_N=\begin{bmatrix}\Phi_N^{(1)} & P\mathcal{L}\\
* & -2\left(\theta_{N+1}^{(1)}-\frac{\lambda_{N+1}}{2\alpha} \right)  \end{bmatrix}<0,\\
&\Phi_N^{(1)} = PF+F^TP+2\delta P+\frac{2\alpha}{\pi^2N}\tilde{K}_0^T\tilde{K}_0.
\end{aligned}
\end{equation}
By Schur complement \eqref{eq:WLMIsNonDelayed} holds iff
\begin{equation}\label{eq:WTailLMINonDelayed}
\scriptsize\begin{bmatrix}\Phi_N^{(1)} & P\mathcal{L} & 0\\
* & -2\theta_{N+1}^{(1)} & 1\\
* & * & -\frac{\alpha}{\lambda_{N+1}} \end{bmatrix}<0.\normalsize
\end{equation}
Note that  LMI \eqref{eq:WTailLMINonDelayed} has $N$-dependent coefficients and dimension.
 Summarizing, we arrive at:
\begin{theorem}\label{Thm:WdynExtension}
Consider  \eqref{eq:PDE1PointActConstChangeVars} with in-domain measurement \eqref{eq:InDomPointMeas1}, control law \eqref{eq:WContDef} and $w(\cdot,0)\in \mathcal{D}(\mathcal{A})$.  Let $\delta>0$ be a desired decay rate, $N_0\in \mathbb{N}$ satisfy \eqref{eq:N0} and $N\in \mathbb{N}$ satisfy $N_0\leq N$. Let $L_0$ and $K_0$ be obtained using \eqref{eq:GainsDesignL}  and \eqref{eq:GainsDesignK}, respectively. Let there exist a $0<P\in \mathbb{R}^{(2N+1)\times (2N+1)}$ and scalar $\alpha>0$ which satisfy \eqref{eq:WTailLMINonDelayed}. Then the solution $w(x,t)$ and $u(t)$ to  \eqref{eq:PDE1PointActConstChangeVars} under the control law \eqref{eq:WContDef}, \eqref{eq:WobsODENonDelayed}
and the corresponding observer $\hat{w}(x,t)$ defined by \eqref{eq:WhatSeries} satisfy
\begin{equation}\label{eq:WH1Stability}
\begin{array}{ll}
&\left\|w(\cdot,t)\right\|^2_{H^1}+ \left|u(t) \right|^2\leq Me^{-2\delta t}\left\|w(\cdot,0)\right\|^2_{H^1},\\[0.5mm]
&\left\|w(\cdot,t)-\hat{w}(\cdot,t)\right\|^2_{H^1}\leq Me^{-2\delta t}\left\|w(\cdot,0)\right\|^2_{H^1}
\end{array}
\end{equation}
with some constant $M>0$. Moreover, \eqref{eq:WTailLMINonDelayed} is always feasible for large enough $N$.
\end{theorem}
\begin{IEEEproof}
Feasibility of \eqref{eq:WTailLMINonDelayed} implies, by the comparison principle,
\begin{equation}\label{eq:ComparisonDirichlet}
V(t)\leq e^{-2\delta t}V(0), \ t\geq 0.
\end{equation}
Since $u(0)=0$, for some constant $M_0>0$ we have
\begin{equation}\label{eq:VZeroNonDelayed1}
\begin{array}{l}
V(0) \overset{\eqref{lem22}}{\leq} M_0\left\|w_x(\cdot,0) \right\|^2_{L^2}\leq M_0\left\|w(\cdot,0)\right\|^2_{H^1}.
\end{array}
\end{equation}
By Wirtinger's inequality (see \cite{Fridman14_TDS}, Sec. 3.10 ),
for $t\geq 0$
\begin{equation}\label{eq:Wirtinger}
\left\|w_x(\cdot,t) \right\|^2\leq \left\|w(\cdot,t) \right\|^2_{H^1}\leq \frac{4+\pi^2}{\pi^2}\left\|w_x(\cdot,t) \right\|^2.
\end{equation}
Since $w(\cdot,t)\in \mathcal{D}(\mathcal{A})$ for all $t> 0$, by \eqref{lem22}
\begin{equation*}
\left\|w_x(\cdot,t) \right\|^2=\sum_{n=1}^{\infty} \lambda_nw_n^2(t).
\end{equation*}
Parseval's equality, \eqref{eq:Wirtinger} and monotonicity of $\lambda_n, \ n\in \mathbb{N}$ imply
\begin{equation}\label{eq:VBoundBelowDirichlet}
\begin{array}{ll}
&V(t)\geq \sigma_{min}(P)\left|u(t) \right|^2+ \frac{\sigma_{min}(P)}{2}\sum_{n=1}^Nw^2_n(t)\\[0.5mm]
&+\sum_{n=N+1}^{\infty}\lambda_nw^2_n(t)\geq \sigma_{min}(P)\left|u(t) \right|^2\\[0.5mm]
&+ \operatorname{min}\left( \frac{\sigma_{min}(P)\pi^2}{2\lambda_N},\pi^2 \right)\left\|w(\cdot,t) \right\|^2_{H^1}, \ \ t\geq 0.
\end{array}
\end{equation}
Then \eqref{eq:WH1Stability} follow from \eqref{eq:ComparisonDirichlet}, \eqref{eq:VZeroNonDelayed1}, \eqref{eq:VBoundBelowDirichlet} and the representation
\begin{equation*}
w(\cdot,t)-\hat{w}(\cdot,t) = \sum_{n=1}^Ne_n(t)\phi_n^D(\cdot)+\sum_{n=N+1}^{\infty}w_n(t)\phi_n^D(\cdot).
\end{equation*}
We will demonstrate feasibility of the derived LMIs for large enough $N$ in the more general setting of $L^2$-gain analysis below (see proof of Theorem \ref{Thm:WdynExtensionHinfDir}). The feasibility of \eqref{eq:WTailLMINonDelayed} for large enough $N$ follows similar arguments.
\end{IEEEproof}
\begin{corollary}\label{cor:H1Conv}
Under the conditions of Theorem \ref{Thm:WdynExtension}, the following estimates hold for $z(x,t)$ satisfying \eqref{eq:ChangeVars}:
\begin{equation}\label{eq:ZH1Stability}
\begin{array}{ll}
&\left\|z(\cdot,t)\right\|^2_{H^1}\leq Me^{-2\delta t}\left\|z(\cdot,0)\right\|^2_{H^1},\\
&\left\|z(\cdot,t)-\hat{w}(\cdot,t)\right\|^2_{H^1}\leq Me^{-2\delta t}\left\|z(\cdot,0)\right\|^2_{H^1},
\end{array}
\end{equation}
where $M>0$ is some constant.
\end{corollary}
\begin{IEEEproof}
From \eqref{eq:ChangeVars} we have
\begin{equation}\label{eq:ZEstimateH1}
\begin{array}{lll}
&\left\|z(\cdot,t) \right\|_{H^1}\leq \left[1+\left\|r(\cdot) \right\|_{H^1}\right]  \max \left(\left\|w(\cdot,t) \right\|_{H^1}, \left|u(t)\right|\right)\\
&\left\|z(\cdot,t) -\hat{w}(\cdot,t)\right\|_{H^1}\leq \left\|w(\cdot,t) -\hat{w}(\cdot,t)\right\|_{H^1}\\
&\hspace{40mm}+\left|u(t)\right|\left\|r(\cdot) \right\|_{H^1}.
\end{array}
\end{equation}
By \eqref{eq:ChangeVars}, \eqref{eq:WH1Stability}, \eqref{eq:ZEstimateH1} and $u(0)=0$ we obtain \eqref{eq:ZH1Stability}.
\end{IEEEproof}
\begin{remark}
Boundary control of 1D heat equation via modal decomposition, without dynamic extension, was considered in \cite{karafyllis2018sampled,katz2020boundary,RamiContructiveFiniteDim}. Without dynamic extension, modal decomposition of KSE \eqref{eq:LinearizedKSE} under boundary conditions \eqref{eq:BCSD} 
results in ODEs similar to \eqref{eq:WOdesNonDelayed} with $v(t)$ replaced by $u(t)$ and $|b_n|\approx \lambda_n^{\frac{3}{2}}$. The growth of $\left\{b_n\right\}_{n=1}^{\infty}$ poses a problem in compensating cross terms (cf. \eqref{eq:WCrosTermNonDelayed}) arising in the Lyapunov stability analysis.  As it is well-known (see e.g. \cite{coron2004global,cerpa2010null}), the use of dynamic extension leads to $\left\{b_n\right\}_{n=1}^{\infty}\in l^2(\mathbb{N})$ (see \eqref{eq:AssbnNonDelayed1}). Similarly, dynamic extension allows to manage with stability analysis under point measurement for the KSE with boundary conditions  \eqref{eq:BCSNe}, and leads to an equivalent control problem with unbounded observation and bounded control operators.
\end{remark}

Let $\delta>0$ and gains $L_0$ and $K_0$ be fixed. The next proposition shows that the feasibility of \eqref{eq:WLMIsNonDelayed} with some $N\geq N_0$ implies the feasibility of \eqref{eq:WLMIsNonDelayed} with $N+1$. In particular, increasing the observer dimension can never result in loss of feasibility (and the decay rate of the closed-loop system for $N$, guaranteed by the LMIs, cannot be better than the one for $N+1$).
\begin{prop}\label{eq:NtoNPlus1}
Let $\delta>0$, $N_0\in \mathbb{N}$ satisfy \eqref{eq:N0} and $N\in \mathbb{N}$ satisfy $N_0\leq N$. Let the gains $L_0$ and $K_0$ be obtained using \eqref{eq:GainsDesignL} and \eqref{eq:GainsDesignK}.
Assume that for some $0<P\in \mathbb{R}^{(2N+1)\times (2N+1)}$ and scalar $\alpha>0$ \eqref{eq:WLMIsNonDelayed} holds with $\theta_{N+1}^{(1)}$ given in \eqref{eq:TailEstDir}. Then, there exists some $0<P_1\in \mathbb{R}^{(2N+3)\times (2N+3)}$ such that \eqref{eq:WLMIsNonDelayed} holds \emph{with $N$ and $P$ replaced by $N+1$ and $P_1$}, respectively, and the same $\alpha>0$.
\end{prop}
\begin{IEEEproof}
Recall $\hat{w}^{N_0}(t)$, $e^{N_0}(t)$, $\hat{w}^{N-N_0}(t)$, $e^{N-N_0}(t)$ and $X_N(t)$ defined in \eqref{eq:WContDef} and \eqref{eq:ErrDefNonDelayed}. For $N+1$, we rewrite $X_{N+1}(t)$ as $X_N(t)$ with the remaining
$e_{N+1}(t),
\hat{w}_{N+1}(t)$ written in the end as follows:
\begin{equation}\label{eq:ChangeofVarsClosedLoop}
\begin{array}{lll}
&\hspace{-2mm}\scriptsize\begin{bmatrix}
X_N(t)\\
e_{N+1}(t)\\
\hat{w}_{N+1}(t)
\end{bmatrix} = \normalsize Q_1 X_{N+1}(t), \  \scriptsize\begin{bmatrix}X_N(t)\\ \zeta_{N+1}(t)\\ e_{N+1}(t) \\ \hat{w}_{N+1}(t) \end{bmatrix} = \normalsize Q_2 \scriptsize\begin{bmatrix}
X_{N+1}(t)\\ \zeta_{N+1}(t)
\end{bmatrix}.
\end{array}
\end{equation}
Here $Q_1$ and $Q_2$ are the following permutation matrices:
\begin{equation}\label{eq:Qmat}
\begin{array}{lll}
& Q_1 = \operatorname{diag}\left\{I_{N+N_0+1},\ \tilde{Q}_1\right\},\ \tilde{Q}_1 = \scriptsize\begin{bmatrix}
0& I_{N-N_0} & 0\\ 0& 0& 1\\ 1&0 &0
\end{bmatrix},\\
&Q_2 = \operatorname{diag}\left\{I_{N+N_0+1},\ \tilde{Q}_2\right\},\ \tilde{Q}_2 = \scriptsize\begin{bmatrix}0 & I_{N-N_0} & 0& 0\\
0 & 0 & 0 & 1\\
0 & 0 & 1 & 0\\
1& 0 & 0 & 0 \end{bmatrix}.
\end{array}
\end{equation}

Let $P_1 = Q_1^T\operatorname{diag} \left\{P,q_1, q_2\right\}Q_1$, where $q_1,q_2>0$ are scalars. Substitute $N+1$ and $P_1$ for $N$ and $P$ in \eqref{eq:VNonDelayed}, respectively. Taking into account the transformations \eqref{eq:ChangeofVarsClosedLoop} and applying arguments similar to \eqref{eq:WStabAnalysisNonDelayed}-\eqref{eq:WStabResultNonDelayed} it can be verified that
\begin{equation}\label{eq:NtoN+1Stab}
\begin{array}{lll}
&\hspace{-2mm}\Psi_{N+1}^{(1)} = \scriptsize Q_2^T \begin{bmatrix}
\Phi_{N+1}^{(1)}  & P\mathcal{L}  & P\mathcal{L}c_{N+1}  & q_2b_{N+1}\tilde{K}_0^T\\
*  & -2\left(\theta_{N+2}^{(1)}-\frac{\lambda_{N+2}}{2\alpha}\right) & 0 & 0 \\
*  & *  &-2q_1\theta_{N+1}^{(1)} & 0\\
*  & *  & *  & -2q_2\theta_{N+1}^{(1)}
\end{bmatrix}Q_2.
\end{array}
\end{equation}
From \eqref{eq:WOdesNonDelayed} we have $b_{N+1}^2={2\over\lambda_{N+1}}$.
 Applying further Schur complement, we find that $\Psi_{N+1}^{(1)}<0$ holds if and only if
\begin{equation*}
\begin{array}{lll}
&S^{(1)}=\scriptsize\begin{bmatrix}
\Phi_{N+1}^{(1)} \ & P\mathcal{L}\\
* \ & -2\left(\theta_{N+2}^{(1)}-\frac{\lambda_{N+2}}{2\alpha}\right)
\end{bmatrix}+\frac{c_{N+1}^2}{2\theta_{N+1}^{(1)}q_1}\scriptsize \begin{bmatrix}
P\mathcal{L}\mathcal{L}^TP \ & 0\\ 0\ & 0
\end{bmatrix}\\
&\hspace{20mm}+\frac{q_2}{\theta_{N+1}^{(1)}\pi^2(N+1)^2}\scriptsize\begin{bmatrix}
\tilde{K}_0^T\tilde{K_0}\ & 0 \\
0 \ & 0
\end{bmatrix}<0.
\end{array}
\end{equation*}
By taking $q_2= 2\alpha\theta_{N+1}^{(1)}$ and $q_1$ sufficiently large we obtain that
$S\leq \Psi_{N}^{(1)}<0.$
\end{IEEEproof}
\begin{remark}\label{rem_bc}
Consider  \eqref{eq:LinearizedKSE} under the different from \eqref{eq:BCSD} and \eqref{eq:BCSNe} boundary conditions
\begin{equation}\label{eq:PDENonDelayedNonLocX2}
\begin{array}{lll}
&z(0,t) =  u(t),  \quad z(1,t)= 0 ,\\
&z_{x}(0,t)=0, \quad z_{x}(1,t)=0.
\end{array}
\end{equation}
Here the eigenfunctions induced by \eqref{eq:SL} are no longer suitable for modal decomposition as their use introduces non-homogeneous terms of the form $z_{xx}(0,t)$ and $z_{xx}(1,t)$ into the ODEs of the modes.
To deal with this difficulty it is theoretically possible to use our approach with the eigenvalues $\left\{\sigma_n\right\}_{n=1}^{\infty}$ and the eigenfunctions $\left\{\psi_n\right\}_{n=1}^{\infty}$ induced by the differential operator $-\partial_{xxxx}-\nu\partial_{xx}$ (see e.g \cite{cerpa2010null,cerpa2017control,guzman2019stabilization}). However, in this case, the eigenvalues $\left\{\sigma_n\right\}_{n=1}^{\infty}$ have neither closed formulas nor estimates of the form (2.3) in \cite{RamiContructiveFiniteDim}.  Instead, they are given as implicit solutions of nonlinear equations, which allow to derive only asymptotic estimates as $n\to \infty$. Similarly, there are no closed formulas for $\left\{\psi_n\right\}_{n=1}^{\infty}$. Note that without closed formulas, the corresponding projections $\bar b_n=<r,\psi_n>$, with $r(x)=1-x$ cannot be computed analytically. It is also not possible to express bounds of the form \eqref{eq:AssbnNonDelayed1} with $b_n$ substituted by $\bar b_n$ (i.e. bounds on $\left\|r\right\|^2-\sum_{n=1}^{N}\left|\left<r, \psi_n\right> \right|^2=\sum_{n=N+1}^{\infty}\left|\left<r, \psi_n\right> \right|^2$) explicitly in terms of $N$.  Hence, for practical implementation, the upper bound in \eqref{eq:AssbnNonDelayed1} can be replaced by the constant $\left\|r\right\|^2$, which may lead to a conservative value of $N$, obtained from LMIs. Moreover, to verify the LMIs feasibility one has to approximate $\left\{<r,\psi_n>\right\}_{n=1}^N$ and $\left\{\sigma_n\right\}_{n=1}^{N+1}$. This large number of numerical approximations can result in a computationally expensive approach with essential numerical errors.
\end{remark}

\emph{B. Neumann actuation and collocated measurement}\label{Sec:2}

Consider the KSE \eqref{eq:LinearizedKSE} with  Neumann boundary conditions \eqref{eq:BCSNe} and
collocated boundary measurement
\begin{equation}\label{eq:BoundMeas}
y(t) = z(0,t).
\end{equation}
Introduce the change of variables
\begin{equation}\label{eq:ChangeVarsNeumann}
w(x,t)=z(x,t)-r(x)u(t), \quad r(x):=x-\frac{x^2}{2}
\end{equation}
to obtain the equivalent
ODE-PDE system
\begin{equation}\label{eq:PDE1PointActConstChangeVarsNeumann}
\begin{aligned}
& \dot{u}(t)=v(t)\\
& w_t(x,t)=-w_{xxxx}(x,t)-\nu w_{xx}(x,t)+\nu u(t)-r(x)v(t),
\end{aligned}
\end{equation}
with 
boundary conditions
\begin{equation}\label{eq:BCsNeumChangeVars}
\begin{array}{lll}
& w_x(0,t)=0, \quad w_x(1,t)=0,\\
& w_{xxx}(0,t) = 0, \quad w_{xxx}(1,t)=0.
\end{array}
\end{equation}
and boundary measurement
\begin{equation}\label{eq:BoundMeas1}
	y(t) = w(0,t).
\end{equation}
Recall that we treat $u(t)$ as an additional state variable and $v(t)$ as the control input, where we choose $u(0)=0$. We present the solution to \eqref{eq:PDE1PointActConstChangeVarsNeumann} as
\begin{equation}\label{eq:WseriesNeumann}
\begin{array}{lll}
w(x,t) &= \sum_{n=0}^{\infty}w_n(t)\phi_n^{Ne}(x), \ w_n(t) &=\left<w(\cdot,t),\phi_n^{Ne}\right>
\end{array}
\end{equation}
with $\left\{\phi_n^{Ne}\right\}_{n\in \mathbb{Z}_+}$ defined in \eqref{eq:SLBCs}. Differentiating under the integral sign, integrating by parts and using \eqref{eq:SL}, \eqref{eq:2BCs} we have
\begin{equation}\label{eq:WOdesNonDelayedNeumann}
\begin{array}{lll}
& \dot{w}_0(t) = \nu u(t)+b_0v(t), \ \ b_0=-\frac{1}{3},\\
&\dot{w}_n(t) = (-\lambda_n^2+\nu\lambda_n)w_n(t) + b_nv(t), \ n\in \mathbb{N} \\
& b_n=-\int_0^1 r(x)\phi^{Ne}_n(x)dx= \frac{\sqrt{2}}{\lambda_n} , \ n\in \mathbb{N},\\
&w_n(0) = \left<w(\cdot,0),\phi_n^{Ne} \right>, \ n\in \mathbb{Z}_+.
\end{array}
\end{equation}
In particular, note that \eqref{eq:AssbnNonDelayed} holds.
Moreover,
\begin{equation}\label{eq:KSENeumannbn}
\sum_{n=N+1}^{\infty}b_n^2\leq \frac{2}{\pi^4}\int_N^{\infty}\frac{1}{x^4}dx\leq \frac{2}{3\pi^4N^3}.
\end{equation}
The faster decay of $b_n, \ n\in \mathbb{Z}_+$, when compared to \eqref{eq:AssbnNonDelayed1}, allows to prove $H^2$-stability of the closed-loop system.

Let $\delta>0$ be a desired decay rate, $N_0 \in \mathbb{Z}_+$ satisfy \eqref{eq:N0} and $N\in \mathbb{Z}_+, \ N_0\leq N$.
We construct a finite-dimensional observer of the form
\begin{equation}\label{eq:WhatSeriesNeumann}
\hat{w}(x,t): = \sum_{n=0}^{N}\hat{w}_n(t)\phi_n^{Ne}(x),
\end{equation}
where $\hat{w}_n(t)$ satisfy the ODEs
\begin{equation}\label{eq:WobsODENonDelayedNeumann}
\begin{array}{lll}
&\dot{\hat{w}}_0(t) = \nu u(t)+b_0v(t) - l_0\left[\hat{w}(0,t)-y(t) \right],\\
&\dot{\hat{w}}_n(t) = (-\lambda_n^2+\nu \lambda_n)\hat{w}_n(t) + b_nv(t)\\
&\hspace{13mm}-l_n\left[\hat{w}(0,t)- y(t)\right],\ n\in \mathbb{N},\\
&\hat{w}_n(0)=0, \quad 0\leq n\leq N.
\end{array}
\end{equation}
with $y(t)$ defined in \eqref{eq:BoundMeas1} and scalar observer gains $\left\{l_n\right\}_{n=0}^N$.\\[0.1cm]
Recall $A_0$ and $\tilde{A}_0$ defined in \eqref{eq:C0A0} and denote
\begin{equation}\label{eq:C0A0Neumann}
\begin{aligned}
&\tilde{A}^{(1)}_0 =\operatorname{diag}\left(\begin{bmatrix}0 & 0\\ \nu & 0 \end{bmatrix} ,A_0\right)\in \mathbb{R}^{(N_0+2)\times(N_0+2)},\\
&L_0^{(1)} = \left[l_0,\dots,l_{N_0} \right]^T,\ \tilde{L}_0^{(1)} = \text{col}\left\{0,L_0^{(1)} \right\}\in \mathbb{R}^{N_0+2},\\
&C_0^{(1)}=\left[c_0,\dots,c_{N_0} \right],\ \tilde{B}_0^{(1)}= \left[1,b_0,\dots,b_{N_0} \right]^T,\\
&c_0 = 1, \ \ c_n= \phi_n^{Ne}(0)=\sqrt{2}, \  n\geq 1.
\end{aligned}
\end{equation}
Assumption 2 and $c_n\neq 0, \ n\geq 0$ imply that the pair $(\tilde{A}_0,C_0^{(1)})$ is observable by the Hautus lemma. Let $l_0,\dots, l_{N_0}$ be such that $L_0^{(1)}$ satisfies the following Lyapunov inequality \begin{equation}\label{eq:GainsDesignLNeumann}
P_{\text{o}}(\tilde{A}_0-L_0^{(1)}C_0^{(1)})+(\tilde{A}_0-L_0^{(1)}C_0^{(1)})^TP_{\text{o}} < -2\delta P_{\text{o}}
\end{equation}
with $0<P_{\text{o}}\in \mathbb{R}^{(N_0+1)\times (N_0+1)}$. Let $l_n=0, \ n>N_0$.

Assumption 2 and \eqref{eq:AssbnNonDelayed} imply that the pair $(\tilde{A}_0^{(1)},\tilde{B}_0^{(1)})$ is controllable. Let $K_0\in \mathbb{R}^{1\times (N_0+2)}$ satisfy
\begin{equation}\label{eq:GainsDesignKNeumann}
\begin{aligned}
&P_{\text{c}}(\tilde{A}_0^{(1)}+\tilde{B}_0^{(1)}K_0)+(\tilde{A}_0^{(1)}+\tilde{B}_0^{(1)}K_0)^TP_{\text{c}} < -2\delta P_{\text{c}},
\end{aligned}
\end{equation}
with $0<P_{\text{c}}\in \mathbb{R}^{(N_0+2)\times (N_0+2)}$.

We propose a $(N_0+2)$-dimensional controller of the form
\begin{equation}\label{eq:WContDefNeumann}
\begin{aligned}
& v(t)= K_0\hat{w}^{N_0}(t),\ \hat{w}^{N_0}(t) = \left[u(t),\hat{w}_0(t),\dots,\hat{w}_{N_0}(t) \right]^T,
\end{aligned}
\end{equation}
which is based on the $N+1$-dimensional observer \eqref{eq:WobsODENonDelayedNeumann}.

For well-posedness of the closed-loop system \eqref{eq:PDE1PointActConstChangeVarsNeumann} and \eqref{eq:WOdesNonDelayedNeumann} with control input \eqref{eq:WContDefNeumann} we consider the operator $\mathcal{A}:\mathcal{D}(\mathcal{A})\to L^2(0,1)$ given in \eqref{eq:DiffOpKSE} and
\begin{equation}\label{eq:DA2}
\mathcal{D}(\mathcal{A})= \left\{h\in H^4(0,1)| h'(0)=h'(1)=h'''(0)=h'''(1)\right\}.
\end{equation}
Integration by parts implies that for $h\in \mathcal{D}(\mathcal{A})$
\begin{equation}\label{eq:CalANeu}
\mathcal{A}h \overset{L^2(0,1)}{=} \sum_{n=1}^{\infty}\left(\lambda_n^2-\nu \lambda_n\right)\left<h,\phi_n^{Ne} \right>\phi_n^{Ne}.
\end{equation}
Let $\mathcal{H}:=L^2(0,1)\times \mathbb{R}^{N+1}$ be a Hilbert space with the norm $\left\|\cdot\right\|_{\mathcal{H}}^2:=\left\| \cdot\right\|^2+\left|\cdot \right|^2$. Introducing the state
\begin{equation*}
	\begin{array}{lll}
		&\xi(t) = \text{col}\left\{\xi^{(1)}(t),\xi^{(2)}(t)\right\},\\
		&\xi^{(1)}(t) = w(\cdot,t), \ \xi^{(2)}(t) = \text{col}\left\{u(t),\hat{w}_0(t),\dots, \hat{w}_N(t) \right\},
	\end{array}
\end{equation*}
the closed-loop system can be presented as \eqref{eq:AbstractODEDirichlet}, \eqref{A1} with $\tilde{A}_0$, $\tilde{B}_0$ and $C_0$ replaced by $\tilde{A}_0^{(1)}$, $\tilde{B}_0^{(1)}$ and $C_0^{(1)}$, respectively. Moreover, $f_2(\xi)$ is now given by $f_2(\xi)= \text{col}\left\{L_0^{(1)}w(0,t),0\right\}$.
The operator $-\tilde{\mathcal{A}}$ generates an analytic semigroup on $\mathcal{H}$. The function $F:\mathcal{D}(\mathcal{A})\times \mathbb{R}^{N+1}\to \mathcal{H}$ is linear. Given $\xi\in \mathcal{D}(\mathcal{A})\times \mathbb{R}^{N+1}$, the Sobolev inequality together with assumption 2 imply
\begin{equation}\label{eq:f2Neumm}
\begin{array}{lll}
&\left\|f_2(\xi) \right\|^2\leq \text{const}_1\cdot\left[\left\|w(\cdot,t) \right\|^2+ \left\|w_{x}(\cdot,t) \right\|^2\right]\\
&\hspace{12mm}\overset{\eqref{lem22}}{=} \text{const}_1\cdot \sum_{n=0}^{\infty}\left[1+\lambda_n\right]\left|\left<w(\cdot,t),\phi_n^{Ne}\right> \right|^2\\
&\hspace{12mm}\overset{\eqref{eq:CalANeu}}{\leq} \text{const}_2\cdot \left[\left\|w(\cdot,t)\right\|^2 + \left\|\mathcal{A}w(\cdot,t)\right\|^2\right]\\
&\hspace{12mm}\leq \text{const}_2 \cdot\left[\left\|\xi\right\|_{\mathcal{H}}^2+\left\|\tilde{\mathcal{A}}\xi\right\|_{\mathcal{H}}^2\right].
\end{array}
\end{equation}
By Theorems 6.3.1 and 6.3.3 in \cite{pazy1983semigroups}, the system \eqref{eq:PDE1PointActConstChangeVarsNeumann}, \eqref{eq:WOdesNonDelayedNeumann} with control input \eqref{eq:WContDefNeumann} and $w(\cdot,0)\in \mathcal{D}(\mathcal{A})$ has a unique classical solution $\xi(t)$, satisfying \eqref{eq:ClassicalSolDirichlet1} and \eqref{eq:ClassicalSolDirichlet2}.

Let \eqref{eq:WEstErrorNonDelayed} be the estimation error. By using \eqref{eq:WseriesNeumann} and \eqref{eq:WhatSeriesNeumann}, the last term on the right-hand side of \eqref{eq:WobsODENonDelayedNeumann} can be written as
\begin{equation}\label{eq:WIntroZetaNonDelayedNeumann}
\begin{array}{ll}
&\hat{w}(0,t)- y(t) =-\sum_{n=0}^{N} c_ne_n(t)-\zeta_N(t),
\end{array}
\end{equation}
where
\begin{equation}\label{eq:zetadefNeumann}
\zeta_N(t)=w(0,t)-\sum_{n=0}^Nw_n(t)\phi_n^{Ne}(0).
\end{equation}
Then the error equations have the form
\begin{equation}\label{eq:WenNeumann}
\begin{array}{ll}
& \dot{e}_0(t) = -l_0\left(\sum_{n=0}^{N} c_ne_n(t)+\zeta_N(t)\right),\\
&\dot e_n(t)=(-\lambda_n^2+\nu \lambda_n)e_n(t)\\
&\hspace{6mm}-l_n\left(\sum_{n=0}^{N} c_ne_n(t)+\zeta_N(t)\right), \ \  1\leq n \leq N_0,\\
& \dot{e}_n(t)=(-\lambda_n^2+\nu \lambda_n)e_n(t), \ \ N_0+1\leq n \leq N.
\end{array}
\end{equation}
To bound $\zeta_N(t)$, given in \eqref{eq:zetadefNeumann}, let $g(x,t):= w(x,t)-\sum_{n=0}^Nw_n(t)\phi_n^{Ne}(x)$ and $\Gamma>0$. By Sobolev's inequality
\begin{equation}\label{eq:zetaEstNeumann1}
\begin{array}{lll}
\zeta_{N}^2(t)\leq &(1+\Gamma)\left\|g(\cdot,t)\right\|^2+ \Gamma^{-1}\left\|g_x(\cdot,t) \right\|^2.
\end{array}
\end{equation}
Using \eqref{eq:SL}, \eqref{eq:2BCs}, \eqref{eq:H2} and \eqref{eq:zetaEstNeumann1} we obtain
\begin{equation}\label{eq:zetaEstNeumann3}
	\zeta_N^2(t)\leq \sum_{n=N+1}^{\infty}\mu_nw_n^2(t), \quad \mu_n=1+\Gamma+\frac{1}{\Gamma }\lambda_n.
\end{equation}
Recall $e^{N-N_0}(t)$, $w^{N-N_0}(t)$, $X_N(t)$, $\tilde{A}_0$, $A_1$, $B_1$ and $C_1$ defined in \eqref{eq:C0A0}, \eqref{A1} and \eqref{eq:ErrDefNonDelayed}. Let
\begin{equation}\label{eq:ErrDefNonDelayedNeumann}
	\begin{array}{lllllll}
		&e^{N_0}(t)=\left[e_0(t),\dots,e_{N_0}(t) \right]^T,\ \tilde{K}_0 = \begin{bmatrix} K_0,&0\end{bmatrix} \in \mathbb{R}^{1\times 2N+3}\\
		& \mathcal{L}^{(1)}= \text{col}\left\{\tilde{L}_0^{(1)},-L_0^{(1)},0\right\}\in \mathbb{R}^{2N+3},\\
		& F^{(1)} = \scriptsize\begin{bmatrix}\tilde{A}_0^{(1)}+\tilde{B}_0^{(1)}K_0 & \tilde{L}_0^{(1)}C_0^{(1)} & 0 &\tilde{L}_0^{(1)}C_1 \\ 0 & \tilde{A}_0-L_0^{(1)}C_0^{(1)} & 0 & -L_0^{(1)}C_1\\ B_1K_0 & 0 & A_1 & 0\\
			0 & 0 & 0 & A_1 \end{bmatrix}.
	\end{array}
\end{equation}
Then, by \eqref{eq:WOdesNonDelayedNeumann}, \eqref{eq:WobsODENonDelayedNeumann}, \eqref{eq:WContDefNeumann}, \eqref{eq:WIntroZetaNonDelayedNeumann}, \eqref{eq:WenNeumann} and \eqref{eq:ErrDefNonDelayedNeumann}, we present the
closed-loop system for $t \geq 0$ as follows:
\begin{equation}\label{eq:ClosedLoopNonDelayedNeumann}
\begin{aligned}
&\dot{X}_N(t) = F^{(1)}X_N(t)+\mathcal{L}^{(1)}\zeta_N(t),\\
& \dot{w}_n(t) = (-\lambda_n^2+\nu\lambda_n)w_n(t) +b_n\tilde{K}_0X_N(t), \ \ n>N.
\end{aligned}
\end{equation}
For $H^2$-stability analysis of \eqref{eq:ClosedLoopNonDelayedNeumann}, define the Lyapunov function
\begin{equation}\label{eq:VNonDelayedNeumann}
V(t)=\left|X_N(t)\right|^2_P+\sum_{n=N+1}^{\infty}\lambda_n^2 w_n^2(t),
\end{equation}
with $0<P\in \mathbb{R}^{(2N+3)\times (2N+3)}$. Differentiation of $V(t)$ along the solution to \eqref{eq:ClosedLoopNonDelayedNeumann} gives
\begin{equation}\label{eq:WStabAnalysisNonDelayedNeumann}
	\begin{array}{lll}
		&\dot{V}+2\delta V = X_N^T(t)\left[PF^{(1)} +\left(F^{(1)}\right)^TP+2\delta P\right]X_N(t)\\
		&+2X_N^T(t)P\mathcal{L}^{(1)}\zeta_N(t)+2\sum_{n=N+1}^{\infty}\lambda_n^2 w_n(t)b_n\tilde{K}_0X_N(t)\\
		&+2\sum_{n=N+1}^{\infty}\left(-\lambda_n^4+\nu \lambda_n^3+\delta \lambda_n^2\right)w_n^2(t)
	\end{array}
\end{equation}
where term-by-term differentiation of the series is justified as in \eqref{eq:WStabAnalysisNonDelayed}. The Cauchy-Schwarz inequality, \eqref{eq:AssbnNonDelayed1} and \eqref{eq:WOdesNonDelayedNeumann} imply
\begin{equation}\label{eq:WCrosTermNonDelayedNeumann}
	\begin{array}{lllll}
		&\sum_{n=N+1}^{\infty} 2\lambda_n^2 w_n(t)b_n\tilde{K}_0X_N(t)\\
		& = 2\sum_{n=N+1}^{\infty}\left[\lambda_n^{\frac{3}{2}} w_n(t) \right]\cdot \left[\frac{\sqrt{2}}{\sqrt{\lambda_n}} \left|\tilde{K}_0X_N(t)\right|\right]\\
		&\leq \frac{1}{\alpha} \sum_{n=N+1}^{\infty}\lambda_n^3 w_n^2(t)+ \frac{2\alpha }{\pi^2N}\left|\tilde{K}_0X_N(t)\right|^2
	\end{array}
\end{equation}
where $\alpha>0$. By monotonicity of $\left\{\lambda_n\right\}_{n=1}^{\infty}$ we have
\begin{equation}\label{eq:TailEstNeum}
\begin{array}{lll}
&2\sum_{n=N+1}^{\infty}\left[-\lambda_n^4+\left(\nu+\frac{1}{2\alpha}\right) \lambda_n^3+\delta\lambda_n^2\right]w_n^2(t)\\
&\leq -2\left(\theta_{N+1}^{(2)}-\frac{\lambda_{N+1}^3}{2\alpha \mu_{N+1}}\right)\sum_{n=N+1}^{\infty}\mu_n w_n^2(t) \\
&\overset{\eqref{eq:zetaEstNeumann3}}{\leq} -2\left(\theta_{N+1}^{(2)}-\frac{\lambda_{N+1}^3}{2\alpha \mu_{N+1}}\right)\zeta_N^2(t),\ \theta_n^{(2)}=\frac{\lambda_n^4-\nu\lambda_n^3-\delta \lambda_n^2}{\mu_n}
\end{array}
\end{equation}
if $-\theta_{N+1}^{(2)}+\frac{\lambda_{N+1}^3}{2\alpha \mu_{N+1}} \leq 0$. Let $\eta(t) = \text{col}\left\{X_N(t),\zeta_N(t) \right\}$. From \eqref{eq:WStabAnalysisNonDelayedNeumann},  \eqref{eq:WCrosTermNonDelayedNeumann} and \eqref{eq:TailEstNeum} we obtain
\begin{equation}\label{eq:WStabResultNonDelayedNeum}
\begin{array}{ll}
&\dot{V}+2\delta V \leq \eta^T(t)\Psi^{(2)} \eta(t)
\end{array}
\end{equation}
provided
\begin{equation}\label{eq:WLMIsNonDelayedNeumann}
\begin{aligned}
&\Psi_N^{(2)}=\begin{bmatrix}\Phi_N^{(2)} & P\mathcal{L}^{(1)}\\
* & -2\left(\theta_{N+1}^{(2)}-\frac{\lambda_{N+1}^3}{2\alpha \mu_{N+1}}\right)  \end{bmatrix}<0,\\
& \Phi_N^{(2)} = PF^{(1)}+\left(F^{(1)}\right)^TP+2\delta P+\frac{2\alpha}{\pi^2N}\tilde{K}_0^T\tilde{K}_0.
\end{aligned}
\end{equation}
Applying Schur complement we have that \eqref{eq:WLMIsNonDelayedNeumann} holds iff
\begin{equation}\label{eq:WTailLMINonDelayedNeumann1}
\begin{bmatrix}\Phi_N^{(2)} & P\mathcal{L}^{(1)} & 0\\
* & -2\theta_{N+1}^{(2)} & 1\\
* & * & -\frac{\alpha \mu_{N+1}}{\lambda_{N+1}^3} \end{bmatrix}<0, \  \mu_{N+1}=1+\Gamma+\frac{1}{\Gamma }\lambda_{N+1}.
\end{equation}
Summarizing, we arrive at:
\begin{theorem}\label{Thm:WdynExtensionNeumann}
Consider  \eqref{eq:PDE1PointActConstChangeVarsNeumann} with  measurement \eqref{eq:BoundMeas1}, control law \eqref{eq:WContDefNeumann} and $w(\cdot,0)\in \mathcal{D}\left(\mathcal{A}\right)$.  Let $\delta>0$ be a desired decay rate, $N_0\in \mathbb{N}$ satisfy \eqref{eq:N0} and $N\in \mathbb{N}$ satisfy $N_0\leq N$. Assume that $L_0$ and $K_0$ are obtained using \eqref{eq:GainsDesignLNeumann}  and \eqref{eq:GainsDesignKNeumann}, respectively. Given $\Gamma>0$, let there exist a positive definite matrix $P\in \mathbb{R}^{(2N+3)\times (2N+3)}$ and scalars $\alpha>0$ such that \eqref{eq:WTailLMINonDelayedNeumann1} holds. Then the solution $w(x,t)$ and $u(t)$ to \eqref{eq:PDE1PointActConstChangeVarsNeumann} under the control law \eqref{eq:WContDefNeumann}, \eqref{eq:WobsODENonDelayedNeumann} and the corresponding observer $\hat{w}(x,t)$ defined by \eqref{eq:WhatSeriesNeumann}
satisfy
\begin{equation}\label{eq:WH2Stability}
\begin{array}{ll}
&\left\|w(\cdot,t)\right\|^2_{H^2}+\left|u(t) \right|^2\leq Me^{-2\delta t}\left\|w(\cdot,0)\right\|^2_{H^2},\\
&\left\|w(\cdot,t)-\hat{w}(\cdot,t)\right\|^2_{H^2}\leq Me^{-2\delta t}\left\|w(\cdot,0)\right\|^2_{H^2}
\end{array}
\end{equation}
with some constant $M>0$. Moreover, \eqref{eq:WTailLMINonDelayedNeumann1} is always feasible for large enough $N$.
\end{theorem}
\begin{IEEEproof}
Feasibility of \eqref{eq:WTailLMINonDelayedNeumann1} implies, by the comparison principle, that \eqref{eq:ComparisonDirichlet} holds. Since $u(0)=0$ and $z(\cdot,0)=w(\cdot,0)\in \mathcal{D}\left(\mathcal{A} \right)$, by Lemma \ref{lem:H1Equiv} we obtain
\begin{equation}\label{eq:VZeroNonDelayedNeumm}
\begin{array}{l}
V(0) \leq M_0 \sum_{n=0}^{\infty}\left(1+\lambda_n^2 \right)w_n^2(t)\overset{\eqref{eq:H2}}{\leq}M_0\left\|z(\cdot,0)\right\|_{H^2}^2
\end{array}
\end{equation}
for some constant $M_0>0$.
Since $w_x(\cdot,t)$ satisfies \eqref{eq:BCsNeumChangeVars}, by Wirtinger's inequality we have
\begin{equation}\label{eq:WirtingerNeumm}
	\left\|w_x(\cdot,t) \right\|^2\leq \frac{4}{\pi^2} \left\|w_{xx}(\cdot,t) \right\|^2.
\end{equation}
  Then, given $t>0$, by arguments similar to \eqref{eq:VBoundBelowDirichlet} we obtain
\begin{equation}\label{eq:VBoundBelowNeumm}
\begin{array}{ll}
&V(t)\geq \sigma_{min}(P)\left|u(t) \right|^2+  M_1\sum_{n=0}^{\infty}\left(1+\lambda_n^2 \right)w^2_n(t)\\[0.5mm]
&\overset{\eqref{eq:H2}}{=} \sigma_{min}(P)\left|u(t) \right|^2+ M_1\left\|w(\cdot,t) \right\|^2+ M_1\left\|w_{xx}(\cdot,t) \right\|^2\\
&\overset{\eqref{eq:WirtingerNeumm}}{\geq} \sigma_{min}(P)\left|u(t) \right|^2+ M_2\left\|w(\cdot,t) \right\|_{H^2}^2
\end{array}
\end{equation}
for some constants $M_1,M_2>0$. 
The rest of the proof follows arguments of Theorem \ref{Thm:WdynExtension}.
\end{IEEEproof}

Similarly to Corollary \ref{cor:H1Conv}, we arrive at
\begin{corollary}
Under the conditions of Theorem \ref{Thm:WdynExtensionNeumann},  the following estimates hold for $z(x,t)$ satisfying \eqref{eq:ChangeVarsNeumann}:
\begin{equation}\label{eq:ZH1StabilityNeumann}
\begin{array}{ll}
&\left\|z(\cdot,t)\right\|^2_{H^2}\leq Me^{-2\delta t}\left\|z(\cdot,0)\right\|^2_{H^2},\\
&\left\|z(\cdot,t)-\hat{w}(\cdot,t)\right\|^2_{H^2}\leq Me^{-2\delta t}\left\|z(\cdot,0)\right\|^2_{H^2}
\end{array}
\end{equation}
with some constant $M>0$.
\end{corollary}
\begin{remark}
By using arguments similar to Proposition \ref{eq:NtoNPlus1} it can be shown that, given $\delta>0$ and gains $L_0$ and $K_0$, the feasibility of \eqref{eq:WTailLMINonDelayedNeumann1} for some $N\geq N_0$ implies the feasibility of \eqref{eq:WTailLMINonDelayedNeumann1} for $N+1$.
\end{remark}

\section{Control with guaranteed $L^2$ and ISS gains}
\emph{A. Dirichlet actuation and in-domain point measurement}\label{Sec:3}\\[0.1cm]
We consider a perturbed version of the PDE \eqref{eq:LinearizedKSE}
\begin{equation}\label{eq:LinearizedKSEHinf}
z_t(x,t) = -z_{xxxx}(x,t) -\nu z_{xx}(x,t)+d(x,t),
\end{equation}
with boundary conditions \eqref{eq:BCSD} and in-domain point measurement
\begin{equation}\label{eq:InDomPointMeasHinf}
y(t) = z(x_*,t)+\sigma(t), \ x_*\in (0,1).
\end{equation}
Here, we consider  disturbances satisfying
\begin{equation}\label{eq:AssumpDisturb}
\begin{array}{lll}
& d\in L^2((0,\infty);L^2(0,1))\cap H^1_{\text{loc}}((0,\infty);L^2(0,1)),\\
& \sigma \in L^2(0,\infty)\cap H^1_{\text{loc}}(0,\infty).
\end{array}
\end{equation}
Introducting the change of variables \eqref{eq:ChangeVars}, we obtain the ODE-PDE system
\begin{equation}\label{eq:PDEDirHinf}
\begin{array}{lll}
& \dot{u}(t)=v(t),\\
&w_t(x,t)=-w_{xxxx}(x,t)-\nu w_{xx}(x,t)-r(x)v(t)+d(x,t)
\end{array}
\end{equation}
with boundary conditions \eqref{eq:BCsDirChangeVars} and measurement
\begin{equation}\label{eq:InDomPointMeas1Hinf}
y(t) = w(x_*,t)+r(x_*)u(t)+\sigma(t).
\end{equation}
Recall that we treat $u(t)$ as an additional state variable and $v(t)$ as the control input, where $u(0)=0$.

We present the solution to \eqref{eq:PDEDirHinf} as \eqref{eq:Wseries}, where $\left\{\phi_n^D\right\}_{n\in \mathbb{N}}$ are defined in \eqref{eq:SLBCs}. Differentiating under the integral sign, integrating by parts and using \eqref{eq:SL} and \eqref{eq:2BCs} we have
\begin{equation}\label{eq:WOdesNonDelayedHInf}
\begin{array}{lll}
\dot{w}_n(t) & = (-\lambda_n^2+\nu\lambda_n)w_n(t) + b_nv(t)+d_n(t), \\
w_n(0)&=\left<w(\cdot,0),\phi_n^D\right>, \   d_n(t)=\left<d(\cdot,t),\phi_n^D\right>
\end{array}
\end{equation}
and $b_n$ defined in \eqref{eq:WOdesNonDelayed} satisfying \eqref{eq:AssbnNonDelayed1}.

We construct a finite-dimensional observer of the form \eqref{eq:WhatSeries}, where $\hat{w}_n(t)$ satisfy \eqref{eq:WobsODENonDelayed} with $y(t)$ defined in \eqref{eq:InDomPointMeas1Hinf} and scalar observer gains $l_n, \ 1\leq n \leq N$.
Let Assumptions 1 and 2 hold. Then the observer and controller gains $L_0$ and $K_0$ can be chosen to satisfy \eqref{eq:GainsDesignL} and \eqref{eq:GainsDesignK}. Let $l_n=0, \ N_0+1\leq n\leq N$.
We propose a $(N_0+1)$-dimensional controller of the form \eqref{eq:WContDef} which is based on the $N$-dimensional observer \eqref{eq:WobsODENonDelayed}.
Then the  closed-loop ODE-PDE system is given by \eqref{eq:PDEDirHinf}, \eqref{eq:BCsDirChangeVars}, \eqref{eq:WobsODENonDelayed} with controller of the form \eqref{eq:WContDef}.
Well-posedness of the closed-loop system \eqref{eq:PDEDirHinf}, \eqref{eq:WobsODENonDelayed} with $y(t)$ defined in \eqref{eq:InDomPointMeas1Hinf} and controller \eqref{eq:WContDef}, under the assumption \eqref{eq:AssumpDisturb} on the  disturbances $d(x,t)$ and $\sigma(t)$ follows by arguments similar to \eqref{eq:DiffOpKSE}-\eqref{eq:ClassicalSolDirichlet2}. Indeed, by assumption \eqref{eq:AssumpDisturb} the non-homogeneous term $F(\xi,t)$ in \eqref{eq:AbstractODEDirichlet} is locally H\"{o}lder continuous and satisfies the condition of Theorem 6.3.3 in \cite{pazy1983semigroups}. Therefore, if $w(\cdot,0)\in \mathcal{D}(\mathcal{A})$ there exists a unique classical solution satisfying \eqref{eq:ClassicalSolDirichlet1} and \eqref{eq:ClassicalSolDirichlet2}.

Let $\gamma>0$ and $\rho_{w},\rho_{u}\geq0$ be scalars. We introduce the performance index
\begin{equation}\label{eq:PerfInd}
\begin{array}{lll}
&J(\rho_w,\rho_u,\gamma)=\int_0^{\infty} \left[\rho_w^2\left\|w(\cdot,t)\right\|_{L^2}^2+\rho_u^2u^2(t)\right.\\
&\hspace{25mm}\left.-\gamma^2\left(\left\|d(\cdot,t)\right\|_{L^2}^2 +\sigma^2(t)\right) \right]dt.
\end{array}
\end{equation}
The closed-loop ODE-PDE system \eqref{eq:PDEDirHinf}, \eqref{eq:BCsDirChangeVars}, \eqref{eq:WobsODENonDelayed}, \eqref{eq:WContDef} has
$L^2$-gain less or equal to $\gamma$ if $J(\rho_w,\rho_u,\gamma)\leq0$ for all disturbances $d(x,t)$ and $\sigma$(t) satisfying \eqref{eq:AssumpDisturb}
along the solutions of the closed-loop system
starting from $w(\cdot,0)\equiv 0$.

We will find conditions that guarantee  that the following inequality holds along the closed-loop system:
\begin{equation}\label{eq:NegVPerfInd}
\begin{array}{lll}
&\dot{V}+2\delta V +W \leq 0,\\
&W = \rho_w^2\left\|w(\cdot,t)\right\|_{L^2}^2+\rho_u^2u^2(t)-\gamma^2\left(\left\|d(\cdot,t)\right\|_{L^2}^2 +\sigma^2(t)\right)
\end{array}
\end{equation}
with $V(t)$ given in \eqref{eq:VNonDelayed} and $\delta=0$.
Indeed, 
integration of \eqref{eq:NegVPerfInd} in $t$ from $0$ to $\infty$ leads to $J(\rho_w,\rho_u,\gamma)\leq0$ for $w(\cdot,0)\equiv 0$.

In the case of $\delta>0$ and $\rho_w=\rho_u=0$,
 \eqref{eq:NegVPerfInd} and the comparison principle imply 
 ISS of the closed-loop system:
\begin{equation}\label{eq:HInfDirISS}
\begin{array}{r}
V(T)\leq e^{-2\delta T}V(0)+\frac{\gamma^2}{2\delta}\sup_{0\leq t \leq T}\left[\left\|d(\cdot,t) \right\|^2_{L^2} +\sigma^2(t)\right]
\\ \forall T>0.
\end{array}
\end{equation}
Note that due to \eqref{eq:VZeroNonDelayed1} and \eqref{eq:VBoundBelowDirichlet}, \eqref{eq:HInfDirISS} yields for some $\overline{M}>\underline{M}>0$ the following inequality:
\begin{equation}\label{ISS}
\begin{array}{ll}
\underline{M}&\left[|u(t)|^2+\left\|w(\cdot,t)\right\|^2_{H^1}\right]
 {\leq}\overline{M} e^{-2\delta T}\left\|w(\cdot,0)\right\|^2_{H^1}\\ & +\frac{\gamma^2}{2\delta}\sup_{0\leq t \leq T}\left[\left\|d(\cdot,t) \right\|^2_{L^2} +\sigma^2(t)\right]
\quad \forall T>0.
\end{array}
\end{equation}
The latter inequality gives the upper bound ${\gamma\over \sqrt{2\delta}}$ on the ISS gain of the closed-loop system.
\begin{remark}
The performance index \eqref{eq:PerfInd}, expressed in terms of $w(x,t)$ and $u(t)$, is considered for simplicity. Note that for a performance index
\begin{equation}\label{eq:Jbar}
\begin{array}{lll}
&\bar{J}(\bar{\rho}_z,\bar{\rho}_u,\gamma)=\int_0^{\infty} \left[\bar{\rho}_z^2\left\|z(\cdot,t)\right\|_{L^2}^2+\bar{\rho}_u^2u^2(t)\right.\\
&\hspace{25mm}\left.-\gamma^2\left(\left\|d(\cdot,t)\right\|_{L^2}^2 +\sigma^2(t)\right) \right]dt,
\end{array}
\end{equation}
where $\gamma>0$ and $\bar{\rho}_z,\bar{\rho}_u\geq0$, the triangle and Cauchy-Schwarz inequalities imply
\begin{equation*}
\bar{J}(\bar{\rho}_z,\bar{\rho}_u,\gamma) \leq J\left(\sqrt{2}\bar{\rho}_z,\sqrt{\frac{2}{3}\bar{\rho}_z^2+\bar{\rho}_u^2},\gamma\right).
\end{equation*}
Thus, \eqref{eq:NegVPerfInd} with  $\delta=0$, $\rho_w= \sqrt{2}\bar{\rho}_z$ and $\rho_u=\sqrt{\frac{2}{3}\bar{\rho}_z^2+\bar{\rho}_u^2}$ implies  $\bar{J}(\bar{\rho}_z,\bar{\rho}_u,\gamma)\leq0$ for $z(\cdot,0)\equiv 0$.
\end{remark}

Using the estimation error \eqref{eq:WEstErrorNonDelayed} and notations \eqref{eq:Wseries} and \eqref{eq:WhatSeries}, the innovation term $\hat{y}(t) - y(t)$ in \eqref{eq:WobsODENonDelayed} can be presented as
\begin{equation}\label{eq:WIntroZetaNonDelayedHinf}
\begin{array}{ll}
&\hat{y}(t) - y(t)=-\sum_{n=1}^{N} c_ne_n(t)-\zeta_N(t)-\sigma(t),
\end{array}
\end{equation}
where $\zeta_N(t)$ appears in \eqref{eq:zetaintegral} and satisfies \eqref{eq:WzetaEst}. Then the error equations have the form
\begin{equation}\label{eq:WenHinf}
\begin{array}{ll}
&\dot e_n(t)=(-\lambda_n^2+\nu \lambda_n)e_n(t)+d_n(t)\\
&\hspace{6mm}-l_n\left(\sum_{n=1}^{N} c_ne_n(t)+\zeta_N(t)+\sigma(t)\right), \ \  1\leq n \leq N_0,\\
& \dot{e}_n(t)=(-\lambda_n^2+\nu \lambda_n)e_n(t)+d_n(t), \ \ N_0+1\leq n \leq N.
\end{array}
\end{equation}

Using \eqref{eq:WobsODENonDelayed}, \eqref{eq:WContDef}, \eqref{eq:ErrDefNonDelayed}, \eqref{eq:ErrDefNonDelayed1}, \eqref{eq:WOdesNonDelayedHInf} and \eqref{eq:WenHinf}, we present the closed-loop system as
\begin{equation}\label{eq:ClosedLoopNonDelayedHinf}
\begin{aligned}
&\dot{X}_N(t) = FX_N(t)+\mathcal{L}\zeta_N(t)+\mathcal{L}\sigma(t)+d^N(t),\quad t\ge 0,\\
& \dot{w}_n(t) = (-\lambda_n^2+\nu\lambda_n)w_n(t) +b_n\tilde{K}_0X_N(t)+d_n(t), \  n>N.
\end{aligned}
\end{equation}
Here
\begin{equation*}
\begin{array}{lll}
&d^N(t) = \text{col} \left(0,d^{N_0}(t),0,d^{N-N_0}(t) \right),\\
&d^{N_0}(t) = \text{col}\left(d_1(t),\dots,d_{N_0}(t) \right),\\
&d^{N-N_0}(t) = \text{col}\left(d_{N_0+1}(t),\dots,d_{N}(t) \right).
\end{array}
\end{equation*}
Recall that we are interested in determining conditions which guarantee \eqref{eq:NegVPerfInd}, with $V(t)$ given in \eqref{eq:VNonDelayed}. By Parseval's equality $W$ can be presented as
\begin{equation}\label{eq:WHinf}
\begin{array}{lll}
&W = \left|X_N(t)\right|_{\Xi}^2+ \rho_w^2\sum_{n=N+1}^{\infty}w_n^2(t)\\
&\hspace{5mm}-\gamma^2\left|d^N(t) \right|^2-\gamma^2\sum_{n=N+1}^{\infty }d_n^2(t)-\gamma^2 \sigma^2(t)
\end{array}
\end{equation}
with $\Xi= \Xi_1^T\Xi_1$ and
\begin{equation}\label{Xi1}
\begin{array}{lll}
\Xi_1 = \scriptsize\begin{bmatrix}\rho_u & 0 & 0 & 0 & 0\\
0 & \rho_w I_{N0} & \rho_w I_{N0} & 0 & 0\\
0& 0 & 0& \rho_w I_{N-N_0}& \rho_w I_{N-N_0} \end{bmatrix}.
\end{array}
\end{equation}

Differentiating $V(t)$ along the solution to \eqref{eq:ClosedLoopNonDelayedHinf} we have
\begin{equation}\label{eq:WStabAnalysisNonDelayedHinfDir}
\begin{array}{lll}
&\dot{V}+2\delta V = X_N^T(t)\left[PF +F^TP+2\delta P\right]X_N(t)\\
&+2X_N^T(t)P\mathcal{L}\left[\zeta_N(t)+\sigma(t)\right]+2X_N^T(t)Pd^N(t)\\
&+2\sum_{n=N+1}^{\infty}(-\lambda_n^3+\nu \lambda_n^2+\delta\lambda_n)w_n^2(t)\\ &+2\sum_{n=N+1}^{\infty}\lambda_n w_n(t)\left[b_n\tilde{K}_0X_N(t)+d_n(t)\right].
\end{array}
\end{equation}
Furthermore, \eqref{eq:AssbnNonDelayed1} and the Cauchy-Schwarz inequality imply
\begin{equation}\label{eq:WCrosTermNonDelayedHinfDir}
\begin{array}{lllll}
&\hspace{-2mm}\sum_{n=N+1}^{\infty} 2\lambda_n w_n(t)\left[b_n\tilde{K}_0X_N(t)+d_n(t)\right]\\
&\overset{\eqref{eq:AssbnNonDelayed1}}{\leq} \frac{2\alpha }{\pi^2N}\left|\tilde{K}_0X_N(t)\right|^2 +\frac{\alpha+\alpha_1}{\alpha \alpha_1} \sum_{n=N+1}^{\infty}\lambda_n^2 w_n^2(t)\\
&+\alpha_1 \sum_{n=N+1}^{\infty}d_n^2(t).
\end{array}
\end{equation}
where $\alpha,\alpha_1>0$. By using \eqref{eq:WHinf}, \eqref{eq:WStabAnalysisNonDelayedHinfDir} and  \eqref{eq:WCrosTermNonDelayedHinfDir} we find
\begin{equation}\label{eq:WStabAnalysisNonDelayedHinfDir1}
\begin{array}{lll}
&\hspace{-1mm}\dot{V}+2\delta V+W \leq X_N^T(t)\left[PF +F^TP+2\delta P+\frac{2\alpha }{\pi^2N}\tilde{K}_0^T\tilde{K}_0\right.\\
&\hspace{-1mm}\left. +\Xi\right]X_N(t)+2X_N^T(t)P\mathcal{L}\left[\zeta_N(t)+\sigma(t)\right]+2X_N^T(t)Pd^N(t)\\
&\hspace{-1mm}-\gamma^2\left[\sigma^2(t)+\left|d^N(t)\right|^2 \right]+\left(\alpha_1-\gamma^2 \right)\sum_{n=N+1}^{\infty}d_n^2(t)\\
&\hspace{-1mm}+2\sum_{n=N+1}^{\infty}\left(-\theta_n^{(3)}+\frac{\lambda_n}{2\alpha}+\frac{\lambda_n}{2\alpha_1}\right)\lambda_nw_n^2(t),
\end{array}
\end{equation}
where
\begin{equation}\label{eq:Theta3}
\theta_n^{(3)} = \lambda_n^2-\nu\lambda_n-\delta-\frac{\rho_w^2}{2\lambda_n}, \quad n>N.
\end{equation}

By 
 monotonicity of $\left\{\lambda_n\right\}_{n=1}^{\infty}$ we have
  $$\begin{array}{l}
 -\theta_n^{(3)}+\frac{\lambda_n}{2\alpha}+\frac{\lambda_n}{2\alpha_1}\\
 \le -\theta_{N+1}^{(3)}+\frac{\lambda_{N+1}}{2\alpha}+\frac{\lambda_{N+1}}{2\alpha_1}\le 0 \quad \forall n\ge N+1, \end{array}$$
implying due to \eqref{eq:WzetaEst}
\begin{equation}\label{eq:ZetaSeriesHInfDir}
\begin{array}{lll}
&2\sum_{n=N+1}^{\infty}(-\theta_n^{(3)}+\frac{\lambda_n}{2\alpha}+\frac{\lambda_n}{2\alpha_1})\lambda_nw_n^2(t)\\
&\hspace{10mm}\leq -2\left(\theta_{N+1}^{(3)}-\frac{\lambda_{N+1}}{2\alpha}-\frac{\lambda_{N+1}}{2\alpha_1}\right)\zeta_N^2(t).
\end{array}
\end{equation}
Let $\eta(t)=\text{col}\left\{X_N(t),\zeta_N(t),d^N(t),\sigma(t)\right\}$ and $\alpha_1 = \gamma^2$. Then, \eqref{eq:WStabAnalysisNonDelayedHinfDir1} and \eqref{eq:ZetaSeriesHInfDir} imply
\begin{equation*}
\dot{V}+2\delta V +W \leq \eta^T(t)\Psi^{(3)}_N\eta(t)\leq 0
\end{equation*}
provided
\begin{equation}\label{eq:LMIsHInfDir}
\begin{array}{lll}
&\Psi_N^{(3)}=\scriptsize\left[
\begin{array}{c|c}
\begin{matrix}
\Phi_N^{(1)}+\Xi \ & P\mathcal{L}\\
* & -2\left(\theta_{N+1}^{(3)}-\frac{\lambda_{N+1}}{2\alpha}-\frac{\lambda_{N+1}}{2\gamma^2}\right)
\end{matrix} & \begin{matrix} P  \ \ & P\mathcal{L} \\ 0 \ \ & 0 \end{matrix} \\
\hline
* \Tstrut & -\gamma^2 I \Tstrut
\end{array}
\right]<0
\end{array}
\end{equation}
where $\Phi_N^{(1)}$ is defined in \eqref{eq:WLMIsNonDelayed}. Applying Schur complement, we find that \eqref{eq:LMIsHInfDir} holds if and only if
\begin{equation}\label{eq:LMIsHInfDir1}
\begin{array}{lll}
&\hspace{-3mm}\scriptsize\left[
\begin{array}{c|c|c|c}
\begin{matrix}
\Phi_N^{(1)} \ & P\mathcal{L}\\
* & -2\theta_{N+1}^{(3)}
\end{matrix} & \begin{matrix}0\qquad &0\\1\qquad &1 \end{matrix}& \begin{matrix} P  \ \ & P\mathcal{L} \\ 0 \ \ & 0 \end{matrix}&\ \begin{matrix}\ \Xi_1^T\\0 \end{matrix}\\
\hline
*\Tstrut& -\operatorname{diag}\left(\frac{\alpha}{\lambda_{N+1}},\frac{\gamma^2}{\lambda_{N+1}} \right)\Tstrut&\ \bigZero[10]\Tstrut&\ \bigZero[10]\Tstrut\\
\hline
* \Tstrut& * \Tstrut& -\gamma^2 I\Tstrut&\ \bigZero[10]\Tstrut\\
\hline * \Tstrut& * \Tstrut& * \Tstrut& -I\Tstrut
\end{array}
\right]<0.
\end{array}
\end{equation}
Note that if \eqref{eq:LMIsHInfDir1} holds for $\delta=0$, then we obtain internal exponential stability of the closed-loop system with a small enough decay rate $\delta_0>0$. Summarizing, we have:
\begin{theorem}\label{Thm:WdynExtensionHinfDir}
Consider the system \eqref{eq:PDEDirHinf} with boundary conditions \eqref{eq:BCSD}, perturbed in-domain measurement \eqref{eq:InDomPointMeas1Hinf} and control law \eqref{eq:WContDef}, \eqref{eq:WobsODENonDelayed}. Here, $d(x,t)$ and $\sigma(t)$ are disturbances satisfying \eqref{eq:AssumpDisturb}. Let $\delta=0$, $N_0\in \mathbb{N}$ satisfy \eqref{eq:N0} and $N\in \mathbb{N}$ satisfy $N_0\leq N$. Let $L_0$ and $K_0$ be obtained using \eqref{eq:GainsDesignL} and \eqref{eq:GainsDesignK}, respectively. Given $\gamma>0$,  let there exist $0<P\in \mathbb{R}^{(2N+1)\times (2N+1)}$ and scalar $\alpha>0$ such that \eqref{eq:LMIsHInfDir1} holds with $\theta_n^{(3)}$ and $\Xi_1$ given by \eqref{eq:Theta3} and \eqref{Xi1} respectively. Then the above 
system 
is internally exponentially stable and satisfies $J(\rho_w,\rho_u,\gamma)\leq0$ for $w(\cdot,0)\equiv 0$. Given $\rho_w,\rho_u>0$, the LMI  \eqref{eq:LMIsHInfDir1} is feasible for $N$ and $\gamma$ large enough.
\end{theorem}
\begin{IEEEproof}
We will show that \eqref{eq:LMIsHInfDir1} is always feasible for large enough $N$ and $\gamma>0$ . We assume, without loss of generality, $\gamma\geq 1$. First, consider $\Xi=\Xi_1^T\Xi_1$ with $\Xi_1$ given in \eqref{Xi1}. Since $\Xi$ is symmetric, the equality
\begin{equation*}
\left|\Xi\right|=\max_{\left|g\right|\leq 1}\left|g^T\Xi g \right|=\max_{\left|g\right|\leq 1}\left|\Xi_1g \right|^2
\end{equation*}
implies $\left|\Xi\right|= \left|\Xi_1\right|^2\leq \max\left(\rho_u^2,2\rho_w^2 \right)$ is independent of $N$. Thus, there exists $0<\mu \in \mathbb{R}$ large enough such that
\begin{equation}\label{eq:11Term}
-\mu I+\Xi<0
\end{equation}
for all $N\in \mathbb{N}$. Next, note that \eqref{eq:AsscnNonDelayed} and \eqref{eq:AssbnNonDelayed1} imply $\left\{c_n\right\}_{n=1}^{\infty}\in l^{\infty}(\mathbb{N})$ and $\left\{b_n\right\}_{n=1}^{\infty}\in l^2(\mathbb{N})$, respectively.  By arguments of Theorem 3.2 in \cite{RamiContructiveFiniteDim}, there exist some $\kappa>0$ and $\Lambda>0$, independent of $N$, such that $\left|e^{\left(F+\delta I\right)t} \right|\leq \Lambda \cdot \sqrt{N}\left(1+t+t^2\right)e^{-\kappa t}$. Therefore, $P\in \mathbb{R}^{(2N+1)\times(2N+1)}$ which solves the Lyapunov equation
\begin{equation}\label{eq:LyapEqKSE}
P(F+\delta I)+(F+\delta I)^T=-\mu I
\end{equation}
satisfies
\begin{equation}\label{eq:PnormKSE}
\left|P\right|\leq \Lambda_1\cdot N,
\end{equation}
where $0<\Lambda_1 \in \mathbb{R}$ is independent of $N$. Substituting  \eqref{eq:LyapEqKSE}, $\lambda_{N+1}=\pi^2\left(N+1\right)^2$ and $\alpha = 1$ into the top left block of \eqref{eq:LMIsHInfDir} we first show
\begin{equation}\label{eq:11Block}
\scriptsize\begin{bmatrix}
-\mu I+\Xi+\frac{2}{\pi^2 N}\tilde{K}_0^T\tilde{K}_0 \ & P\mathcal{L}\\
* & -2\left(\theta_{N+1}^{(3)}-\frac{\lambda_{N+1}}{2}-\frac{\lambda_{N+1}}{2\gamma^2}\right)
\end{bmatrix}\normalsize<0
\end{equation}
holds for large enough $N$. From \eqref{eq:11Term}, $\gamma \geq 1$ and $\lambda_{N+1}\approx \left( N+1\right)^2$, the diagonal blocks are negative provided $N$ is large enough. Applying Schur complement \eqref{eq:11Block} holds iff
\begin{equation}\label{eq:EquivLMIsKSE}
\begin{array}{lll}
&-\mu I+ \Xi+ \frac{2}{\pi^2N}\tilde{K}_0^T\tilde{K}_0\\
&\hspace{20mm}+\frac{1}{2\left(\theta_{N+1}^{(3)}-\frac{\lambda_{N+1}}{2}-\frac{\lambda_{N+1}}{2\gamma^2}\right)}P\mathcal{L}\mathcal{L}^TP<0.
\end{array}
\end{equation}
Note that $\theta_n^{(3)}\approx n^4$ for large $n$, whereas $\left|\tilde{K}_0\right|$ and $\left|\mathcal{L}\right|$ are independent of $N$. Taking into account \eqref{eq:PnormKSE} and increasing $N$ we have that \eqref{eq:11Block} holds for large enough $N$. Finally, consider \eqref{eq:LMIsHInfDir} with $N$ large enough for \eqref{eq:11Block} to hold. Applying Schur complement and choosing $\gamma$ large enough, \eqref{eq:LMIsHInfDir} holds.
\end{IEEEproof}
\begin{remark}\label{rem:IssHInfDir}
Let \eqref{eq:LMIsHInfDir1} hold with $\delta>0$ and $\rho_w = \rho_u=0$ (i.e $\Xi_1=0$). Then the closed-loop system \eqref{eq:PDEDirHinf}, \eqref{eq:BCSD}, \eqref{eq:WContDef}, \eqref{eq:WobsODENonDelayed} is ISS and its solutions satisfy the bounds \eqref{eq:HInfDirISS} and \eqref{ISS}.
\end{remark}
\begin{remark}
One may want to extend Proposition \ref{eq:NtoNPlus1} to the cases of ISS and $L^2$-gain analysis and show that given $\delta>0$ and $\gamma^2$, the feasibility of LMI \eqref{eq:LMIsHInfDir1} with some $N\geq N_0$ implies the feasibility of LMI \eqref{eq:LMIsHInfDir1} with $N+1$. Differently from stabilization, in ISS and $L^2$-gain a coupling of $e_{N+1}(t)$ and $d_{N+1}(t)$ is given in the ODE of $e_{N+1}(t)$ (see \eqref{eq:WenHinf}). Therefore, $e_{N+1}(t)$ is no longer exponentially decaying. Furthermore, coupling of $X_N(t)$ and $e_{N+1}(t)$ is introduced through the innovation term \eqref{eq:WIntroZetaNonDelayedHinf}. Therefore, for ISS and $L^2$-gain , the proof of Proposition \ref{eq:NtoNPlus1} fails to follow through. Indeed, consider the case of ISS (i.e $\rho_w=\rho_u=0$, which implies $\Xi_1=0$ in \eqref{Xi1}). Recall $Q_1$, given in \eqref{eq:Qmat}, which satisfies \eqref{eq:ChangeofVarsClosedLoop}. As in Proposition \ref{eq:NtoNPlus1}, let $P_1 = Q_1^T\operatorname{diag} \left\{P,q_1, q_2\right\}Q_1$, where $q_1,q_2>0$ are scalars. Substitution of $P_1$, $\alpha>0$, $\gamma>0$ and $\rho_u=\rho_w = 0$ into \eqref{eq:LMIsHInfDir} results in the following equivalent LMI for $\eta(t)=\text{col}\left\{X_N(t),\zeta_{N+1}(t),d^N(t),\sigma(t),e_{N+1}(t),\hat{w}_{N+1}(t),d_{N+1}(t)\right\}$
\begin{equation}\label{eq:NtoN+1}
\scriptsize\left[
\begin{array}{c|c|c}
\begin{matrix}
\Phi_{N+1}^{(1)} \ & P\mathcal{L} \\
* & -2\left(\theta_{N+2}^{(3)}-\frac{\lambda_{N+2}}{2\alpha}-\frac{\lambda_{N+2}}{2\gamma^2}\right)
\end{matrix} & \begin{matrix} P  \ \ & P\mathcal{L} \\ 0 \ \ & 0 \end{matrix} & \begin{matrix}
\Pi_1\ & 0\\
0 \ & 0
\end{matrix}  \\
\hline
* \Tstrut & -\gamma^2 I \Tstrut& 0\\
\hline
* \Tstrut & * \Tstrut & \Pi_2
\end{array}
\right]<0, \normalsize
\end{equation}
where
\begin{equation}\label{eq:NtoN+1_1}
\begin{array}{lll}
&\Pi_1 = \scriptsize\begin{bmatrix}
P\mathcal{L}c_{N+1} \ & \ q_2b_{N+1}\tilde{K}_0^T
\end{bmatrix},\\
&\Pi_2 = \scriptsize\begin{bmatrix}
-2q_1(\lambda_{N+1}^2-\nu\lambda_{N+1}) \ &  0 \  &  q_1\\
* \ &  -2q_2(\lambda_{N+1}^2-\nu\lambda_{N+1}) \ & 0\\
* \ &  * \ &  -\gamma^2
\end{bmatrix},\\
&\Phi_{N+1}^{(1)} = PF+F^TP+2\delta P+\frac{2\alpha}{\pi^2(N+1)}\tilde{K}_0^T\tilde{K}_0.
\end{array}
\end{equation}
Note that unlike \eqref{eq:NtoN+1Stab}, where we had $q_1$ and $q_2$ appearing only on the diagonal, $q_1$ appears off-diagonal in $\Pi_2$. To proceed, we need to verify that $\Pi_2<0$. By Schur complement, the latter holds iff
\begin{equation*}
\scriptsize \begin{bmatrix}
-2q_1(\lambda_{N+1}^2-\nu\lambda_{N+1})+\frac{q_1^2}{\gamma^2} \ & 0\\
* \ & -2q_2(\lambda_{N+1}^2-\nu\lambda_{N+1})
\end{bmatrix}<0.\normalsize
\end{equation*}
In particular, $q_1$ must satisfy $q_1<2\gamma^2(\lambda_{N+1}^2-\nu\lambda_{N+1})$.  Using Schur complement we have that \eqref{eq:NtoN+1} holds iff
\begin{equation}\label{eq:S3}
\begin{array}{lll}
&S^{(3)}=\scriptsize\left[
\begin{array}{c|c}
\begin{matrix}
\Phi_{N+1}^{(1)}+\Upsilon & P\mathcal{L}\\
* & -2\left(\theta_{N+2}^{(3)}-\frac{\lambda_{N+2}}{2\alpha}-\frac{\lambda_{N+2}}{2\gamma^2}\right)
\end{matrix} & \begin{matrix} P  \ \ & P\mathcal{L} \\ 0 \ \ & 0 \end{matrix} \\
\hline
* \Tstrut & -\gamma^2 I \Tstrut
\end{array}
\right]<0,\\
&\Upsilon =\scriptsize \frac{c_{N+1}^2}{2}\left[\frac{1}{2q_1\theta_{N+1}^{(3)}}+ \frac{1}{\theta_{N+1}^{(3)}\left(2\theta_{N+1}^{(3)}-q_1 \right)} \right]P\mathcal{L}\mathcal{L}^TP\\
&\hspace{5mm}+\frac{q_2b_{N+1}^2}{\theta_{N+1}^{(3)}}\tilde{K}_0^T\tilde{K}_0.
\end{array}
\end{equation}
Here we can take $q_2\to 0^+$ small enough. However, due to the condition $q_1<2\gamma^2(\lambda_{N+1}^2-\nu\lambda_{N+1})$, it is not possible to take $q_1\to \infty$. Similar restrictions on $q_1$ are obtained for $L^2$-gain analysis. 

{Note that for the case of ISS with $d(x,t)\equiv 0$ we obtain the equivalent LMI \eqref{eq:NtoN+1} with the last column and row removed. Therefore, no restriction on $q_1$ is imposed. Taking $q_2\to 0^+$ small and $q_1\to \infty$ large enough, we have that feasibility of LMI \eqref{eq:LMIsHInfDir1} with some $N$ implies the feasibility of LMI \eqref{eq:LMIsHInfDir1} with $N+1$. For $L^2$-gain, this remains unclear.}
\end{remark}

\emph{B. Neumann actuation and collocated measurement}\label{Sec:4}

Consider the perturbed version of the PDE \eqref{eq:LinearizedKSE}, given by \eqref{eq:LinearizedKSEHinf},
with disturbances $d(x,t)$ and $\sigma(t)$ satisfying \eqref{eq:AssumpDisturb}, boundary conditions \eqref{eq:BCSNe} and collocated boundary measurement
\begin{equation}\label{eq:InDomPointMeasHinfNeu}
y(t) = z(0,t)+\sigma(t).
\end{equation}
By change of variables \eqref{eq:ChangeVarsNeumann}, we obtain the ODE-PDE system
\begin{equation}\label{eq:PDEDirHinfNeu}
\begin{array}{lll}
&\dot{u}(t)=v(t),\\
&w_t(x,t)=-w_{xxxx}(x,t)-\nu w_{xx}(x,t)\\
&\hspace{20mm}+\nu u(t)-r(x)v(t)+d(x,t)
\end{array}
\end{equation}
with boundary conditions \eqref{eq:BCsNeumChangeVars} and measurement
\begin{equation}\label{eq:InDomPointMeas1HinfNeu}
y(t) = w(0,t)+\sigma(t).
\end{equation}
We present the solution to \eqref{eq:PDEDirHinfNeu} as \eqref{eq:WseriesNeumann}, where $\phi_n^{Ne}(x)$ are defined in \eqref{eq:SLBCs}. Differentiating under the integral sign, integrating by parts and using \eqref{eq:SL} and \eqref{eq:2BCs} we have
\begin{equation}\label{eq:WOdesNonDelayedHInfNeumann}
\begin{array}{lll}
& \dot{w}_0(t) = \nu u(t)+b_0v(t)+d_0(t), \\
&\dot{w}_n(t) = (-\lambda_n^2+\nu\lambda_n)w_n(t) + b_nv(t)+d_n(t), \ n\in \mathbb{N} \\
&w_n(0) = \left<w(\cdot,0),\phi_n^{Ne} \right>, \ n\in \mathbb{Z}_+
\end{array}
\end{equation}
and $b_n, \ n\in \mathbb{Z}_+$ defined in \eqref{eq:WOdesNonDelayedNeumann} satisfy \eqref{eq:KSENeumannbn}.

Let $\delta\geq0$, $N_0 \in \mathbb{Z}_+$ satisfy \eqref{eq:N0} and $N\in \mathbb{Z}_+, \ N_0\leq N$. Let scalars $\gamma>0$ and $\rho_{w},\rho_{u}\geq0$. Recall the performance index given by \eqref{eq:PerfInd}. We are interested in finding a control law $v(t)$ which guarantees \eqref{eq:NegVPerfInd}, where $V(t)$ is given by \eqref{eq:VNonDelayed} with $\lambda_n$ replaced by $\lambda_n^{1.25}$ for $n\geq N+1$.
This choice of $V(t)$ is done to avoid further restrictions on the class of admissible disturbances.

We construct a finite-dimensional observer of the form \eqref{eq:WhatSeriesNeumann} where $\hat{w}_n(t)$ satisfy the ODEs \eqref{eq:WobsODENonDelayedNeumann} with $y(t)$ defined in \eqref{eq:InDomPointMeas1HinfNeu} and scalar observer gains $l_n, \ 0\leq n \leq N$.
Let Assumptions 1 and 2 hold. Then, the observer and controller gains $L_0^{(1)}$ and $K_0$ can be chosen to satisfy \eqref{eq:GainsDesignLNeumann} and \eqref{eq:GainsDesignKNeumann}, respectively. Let $l_n=0, \ N_0+1\leq n\leq N$.

We propose a $(N_0+2)$-dimensional controller of the form \eqref{eq:WContDefNeumann} which is based on the $N+1$-dimensional observer \eqref{eq:WobsODENonDelayedNeumann}.\\[0.1cm]
Using the estimation error $e_n(t)=w_n(t)-\hat{w}_n(t),\ 0\leq n\leq N$, \eqref{eq:WseriesNeumann} and \eqref{eq:WhatSeriesNeumann}, the innovation term $\hat{y}(t) - y(t)$ in \eqref{eq:WobsODENonDelayedNeumann} can be presented as \eqref{eq:WIntroZetaNonDelayedHinf} (with summation starting at $n=0$), where $\zeta_N(t)$ appears in \eqref{eq:zetadefNeumann} and satisfies \eqref{eq:zetaEstNeumann3}. Then the error equations have the form
\begin{equation}\label{eq:WenHinfNeum}
\begin{array}{ll}
& \dot{e}_0(t) = -l_0\left(\sum_{n=0}^{N} c_ne_n(t)+\zeta_N(t)+\sigma(t)\right)+d_0(t),\\
&\dot e_n(t)=(-\lambda_n^2+\nu \lambda_n)e_n(t)+d_n(t)\\
&\hspace{6mm}-l_n\left(\sum_{n=0}^{N} c_ne_n(t)+\zeta_N(t)+\sigma(t)\right), \   1\leq n \leq N_0,\\
& \dot{e}_n(t)=(-\lambda_n^2+\nu \lambda_n)e_n(t)+d_n(t), \ \ N_0+1\leq n \leq N.
\end{array}
\end{equation}

Well-posedness of the closed-loop system \eqref{eq:PDEDirHinfNeu} and \eqref{eq:WobsODENonDelayedNeumann} with \eqref{eq:AssumpDisturb}, $y(t)$ in \eqref{eq:InDomPointMeas1HinfNeu} and control law \eqref{eq:WContDefNeumann} follows from arguments similar to \eqref{eq:DA2}-\eqref{eq:f2Neumm}. Thus, for $w(\cdot,0)\in \mathcal{D}(\mathcal{A})$ there exists a unique classical solution satisfying \eqref{eq:ClassicalSolDirichlet1} and \eqref{eq:ClassicalSolDirichlet2}.

Using \eqref{eq:ErrDefNonDelayed}, \eqref{eq:WobsODENonDelayedNeumann}, \eqref{eq:WContDefNeumann},  \eqref{eq:ErrDefNonDelayedNeumann}, \eqref{eq:WOdesNonDelayedHInfNeumann} and \eqref{eq:WenHinfNeum}, we arrive at the closed-loop system
\begin{equation}\label{eq:ClosedLoopNonDelayedHinfNeum}
\begin{aligned}
&\dot{X}_N(t) = F^{(1)}X_N(t)+\mathcal{L}^{(1)}\zeta_N(t)+\mathcal{L}^{(1)}\sigma(t)+d^N(t),\\
& \dot{w}_n(t) = (-\lambda_n^2+\nu\lambda_n)w_n(t) +b_n\tilde{K}_0X_N(t)+d_n(t), \ \ n>N.
\end{aligned}
\end{equation}
We derive conditions which guarantee \eqref{eq:NegVPerfInd}, with $V(t)$ given by \eqref{eq:VNonDelayed}.

Differentiation of $V(t)$ along the solution to \eqref{eq:ClosedLoopNonDelayedHinfNeum} gives
\begin{equation}\label{eq:WStabAnalysisNonDelayedHinfNeum}
\begin{array}{lll}
&\dot{V}+2\delta V = X_N^T(t)\left[PF +F^TP+2\delta P\right]X_N(t)\\
&+2X_N^T(t)P\mathcal{L}\left[\zeta_N(t)+\sigma(t)\right]+2X_N^T(t)Pd^N(t)\\
&+2\sum_{n=N+1}^{\infty}(-\lambda_n^3+\nu \lambda_n^2+\delta\lambda_n)w_n^2(t)\\ &+2\sum_{n=N+1}^{\infty}\lambda_n w_n(t)\left[b_n\tilde{K}_0X_N(t)+d_n(t)\right].
\end{array}
\end{equation}
By the Cauchy-Schwarz inequality and $b_n$ given in \eqref{eq:WOdesNonDelayedNeumann}
\begin{equation}\label{eq:WCrosTermNonDelayedHinfNeum}
\begin{array}{lllll}
&\hspace{-2mm}\sum_{n=N+1}^{\infty} 2\lambda_n w_n(t)\left[b_n\tilde{K}_0X_N(t)+d_n(t)\right]\\
&\overset{}{\leq} \frac{2\alpha }{\pi^2N}\left|\tilde{K}_0X_N(t)\right|^2 +\frac{\alpha+\alpha_1}{\alpha \alpha_1} \sum_{n=N+1}^{\infty}\lambda_n^2 w_n^2(t)\\
&+\alpha_1 \sum_{n=N+1}^{\infty}d_n^2(t).
\end{array}
\end{equation}
with $\alpha,\alpha_1>0$. By \eqref{eq:WHinf}, \eqref{eq:WStabAnalysisNonDelayedHinfNeum} and \eqref{eq:WCrosTermNonDelayedHinfNeum} we find
\begin{equation}\label{eq:WStabAnalysisNonDelayedHinfNeum1}
\begin{array}{lll}
&\hspace{-1mm}\dot{V}+2\delta V+W \leq X_N^T(t)\left[PF +F^TP+2\delta P+\frac{2\alpha }{\pi^2N}\tilde{K}_0^T\tilde{K}_0\right.\\
&\hspace{-1mm}\left. +\Xi\right]X_N(t)+2X_N^T(t)P\mathcal{L}\left[\zeta_N(t)+\sigma(t)\right]+2X_N^T(t)Pd^N(t)\\
&\hspace{-1mm}-\gamma^2\left[\sigma^2(t)+\left|d^N(t)\right|^2 \right]+\left(\alpha_1-\gamma^2 \right)\sum_{n=N+1}^{\infty}d_n^2(t)\\
&\hspace{-1mm}+2\sum_{n=N+1}^{\infty}\left(-\theta_n^{(1)}+\frac{\lambda_n}{2\alpha}+\frac{\lambda_n}{2\alpha_1}\right)\lambda_nw_n^2(t),
\end{array}
\end{equation}
where $\mu_n, \ n>N$ is defined in \eqref{eq:zetaEstNeumann3} and
\begin{equation}\label{eq:theta4}
\theta_n^{(4)} = \frac{\lambda_n^{3}-\nu \lambda_n^{2}-\delta \lambda_n-0.5\rho_w^2}{\mu_n}, \quad n>N.
\end{equation}
By monotonicity of $\lambda_n, \ n\geq 0$ we have
$$\begin{array}{l}-\theta_n^{(4)}+\frac{\lambda_n}{2\alpha \mu_n}+\frac{\lambda_n^{2}}{2\alpha_1\mu_n}\\
\le-\theta_{N+1}^{(4)}+\frac{\lambda_{N+1}}{2\alpha \mu_{N+1}}+\frac{\lambda_{N+1}^{2}}{2\alpha_1\mu_{N+1}}\leq 0 \quad \forall n>N. \end{array}$$
Then, due to \eqref{eq:zetaEstNeumann3} we obtain
\begin{equation}\label{eq:ZetaSeriesHInfNeum}
\begin{array}{lll}
&2\sum_{n=N+1}^{\infty}\left(-\theta_n^{(4)}+\frac{\lambda_n}{2\alpha \mu_n}+\frac{\lambda_n^{2}}{2\alpha_1\mu_n}\right)\mu_nw_n^2(t)\\
&\hspace{10mm}\leq 2\left(-\theta_{N+1}^{(4)}+\frac{\lambda_{N+1}}{2\alpha \mu_{N+1}}+\frac{\lambda_{N+1}^{2}}{2\alpha_1\mu_{N+1}}\right)\zeta_N^2(t).
\end{array}
\end{equation}

Let $\eta(t)=\text{col}\left(X_N(t),\zeta_N(t),d^N(t),\sigma(t)\right)$ and $\alpha_1 = \gamma^2$. Then, \eqref{eq:WStabAnalysisNonDelayedHinfNeum1} and \eqref{eq:ZetaSeriesHInfNeum} imply
\begin{equation*}
\dot{V}+2\delta V +W \leq \eta^T(t)\Psi_N^{(4)}\eta(t)\leq 0
\end{equation*}
provided
\begin{equation}\label{eq:LMIsHInfNeu}
\begin{array}{lll}
&\hspace{-2mm}\Psi_N^{(4)}=\scriptsize\left[
\begin{array}{c|c}
\begin{matrix}
\Phi_N^{(2)}+\Xi \ & P\mathcal{L}^{(1)}\\
* & -2\left(\theta_{N+1}^{(4)}-\frac{\lambda_{N+1}\left(\gamma^2+\alpha \lambda_{N+1} \right)}{2\alpha \gamma^2 \mu_{N+1}}\right)
\end{matrix} & \begin{matrix} P  \ \ & P\mathcal{L}^{(1)} \\ 0 \ \ & 0 \end{matrix} \\
\hline
*\Tstrut& -\gamma^2 I_2\Tstrut
\end{array}
\right]<0,\\
&\hspace{-2mm} \normalsize \mu_{N+1} =1+\Gamma+\frac{1}{\Gamma }\lambda_{N+1},
\end{array}
\end{equation}
where $\Phi_N^{(2)}$ is defined in \eqref{eq:WLMIsNonDelayedNeumann}. By Schur complement \eqref{eq:LMIsHInfNeu} holds if and only if
\begin{equation}\label{eq:LMIsHInfNeu1}
\begin{array}{lll}
&\scriptsize\left[
\begin{array}{c|c|c|c}
\begin{matrix}
\Phi_N^{(2)} \ & P\mathcal{L}^{(1)}\\
* & -2\theta_{N+1}^{(4)}
\end{matrix} & \begin{matrix}0\qquad &0\\1\qquad &1 \end{matrix}& \begin{matrix} P  \ \ & P\mathcal{L}^{(1)} \\ 0 \ \ & 0 \end{matrix}&\ \begin{matrix}\ \Xi_1^T\\0 \end{matrix}\\
\hline
*\Tstrut& -\operatorname{diag}\left(\frac{\alpha\mu_{N+1}}{\lambda_{N+1}},\frac{\gamma^2\mu_{N+1}}{\lambda_{N+1}^{2}} \right)\Tstrut& \bigZero[10]\Tstrut&\ \bigZero[10]\Tstrut\\
\hline
* \Tstrut& * \Tstrut& -\gamma^2 I\Tstrut&\ \bigZero[10]\Tstrut\\
\hline * \Tstrut& * \Tstrut& * \Tstrut& -I\Tstrut
\end{array}
\right]<0,\\
&\normalsize \mu_{N+1} =1+\Gamma+\frac{1}{\Gamma }\lambda_{N+1} .
\end{array}
\end{equation}
Note that if \eqref{eq:LMIsHInfNeu1} holds for $\delta=0$ then we obtain internal exponential stability of the closed-loop system with a small enough decay rate $\delta_0>0$.

Summarizing, we have:
\begin{theorem}\label{Thm:WdynExtensionHinfNeu}
Consider the system \eqref{eq:PDEDirHinfNeu} with boundary conditions \eqref{eq:BCsNeumChangeVars}, boundary measurement \eqref{eq:InDomPointMeas1HinfNeu} and control law \eqref{eq:WContDefNeumann}. Here, $d(x,t)$ and $\sigma(t)$ are disturbances satisfying \eqref{eq:AssumpDisturb}. Let $\delta=0$, $N_0\in \mathbb{N}$ satisfy \eqref{eq:N0} and $N\in \mathbb{N}$ satisfy $N_0\leq N$. Let $L_0$ and $K_0$ be obtained using \eqref{eq:GainsDesignLNeumann} and \eqref{eq:GainsDesignKNeumann}. Given $\gamma>0$ and $\Gamma>0$, let there exist $0<P\in \mathbb{R}^{(2N+2)\times (2N+2)}$ and a scalar $\alpha>0$ satisfying \eqref{eq:LMIsHInfNeu1} with $\theta_n^{(4)}$ given by \eqref{eq:theta4}. Then \eqref{eq:PDEDirHinfNeu} is internally exponentially stable and satisfies $J(\rho_w,\rho_u,\gamma)\leq0$ for $w(\cdot,0)\equiv 0$. Furthermore, given $\rho_w,\rho_u>0$, the LMI  \eqref{eq:LMIsHInfDir1} is always feasible for $N$ and $\gamma>0$ large enough.
\end{theorem}

\section{Examples}
We consider KSE  \eqref{eq:LinearizedKSE} with $\nu = 10$. This choice corresponds to an unstable open-loop system in both cases of Dirichlet and Neumann actuations.
The feasibility of LMIs is verified using the Matlab LMI toolbox. 
\\[0.1cm]
\emph{A. Dirichlet actuation and in-domain measurement}\label{ex:KSEDir}

Consider the unperturbed KSE \eqref{eq:LinearizedKSE},  boundary conditions \eqref{eq:BCSD} and measurement \eqref{eq:InDomPointMeas} with $x_*=\pi^{-1}$. Let $\delta=1$, which results in $N_0=1$. {To guarantee a small value of $N$, the gains $K_0$ and $L_0$ were found by solving \eqref{eq:GainsDesignL} and \eqref{eq:GainsDesignK} with strong inequality replaced by equality and $\delta=\delta_0 \in (1,10)$. The minimal value of $N$ was obtained for $\delta_0 = 5$. The corresponding gains are given by}
\begin{equation}\label{eq:GainsDirStab}
K_0 = \left[7.1415, 26.0901 \right], \ L_0=2.3419.
\end{equation}
The LMI of Theorem \ref{Thm:WdynExtension} is feasible for minimal $N=4$.

Next, we consider the perturbed KSE \eqref{eq:LinearizedKSEHinf} under boundary conditions \eqref{eq:BCSD} and preturbed measurement \eqref{eq:InDomPointMeasHinf} with $x_*=\pi^{-1}$. Here, the disturbances $d(x,t)$ and $\sigma(t)$ satisfy \eqref{eq:AssumpDisturb}.
For the case of input-to-state stabilization we choose $K_0$ and $L_0$ given by \eqref{eq:GainsDirStab}.
For the corresponding $L^2$-gain problem we consider $\rho_w= 0.1$, $\rho_u = 0.2$ and $\delta=0$. {Similarly to \eqref{eq:GainsDirStab}, the gains $K_0$ and $L_0$ were found by solving \eqref{eq:GainsDesignL} and \eqref{eq:GainsDesignK} with strong inequality replaced by equality and $\delta=\delta_0=1.5$. The resulting gains are given by}
\begin{equation}\label{GainsL2}
K_0 = \left[3.0672, 15.911 \right], \ L_0=1.501.
\end{equation}
The LMI \eqref{eq:LMIsHInfDir1}  (with $\delta=0$ and gains \eqref{GainsL2} for $L^2$-gain analysis and with $\delta=1$ and gains \eqref{eq:GainsDirStab} for ISS) is verified for $N \in \left\{4,6,8,10,12\right\}$. For each choice of $N$, we find the smallest $\gamma$ which guarantees the feasibility of the LMI. The results are presented in Table \ref{Tab:SimISSDir}. Note that for ISS $\gamma$ decreases as $N$ grows, whereas for $L^2$-gain the resulting $\gamma$ does not grow for larger $N$.
\begin{table}
	\begin{center}
		\scalebox{0.70}{
			\begin{tabular}{|l|l|l|l|l|l|l|l|}
				\hline
				N & 4 & 6 & 8 &10 & 12 \\
				\hline
				$\gamma$ (ISS) &0.8 & 0.5 & 0.3 & 0.3 & 0.2 \\
				\hline
				$\gamma$ ($L^2\text{-gain}$) &15 & 15 & 15 & 15 & 15 \\
				\hline
		\end{tabular}}
	\end{center}
	\caption{\label{Tab:SimISSDir} Feasibility of LMIs - Dirichlet. $N$ vs minimal $\gamma$. }
\end{table}

\begin{figure}
	\centering
	\includegraphics[width=48mm,scale=0.11]{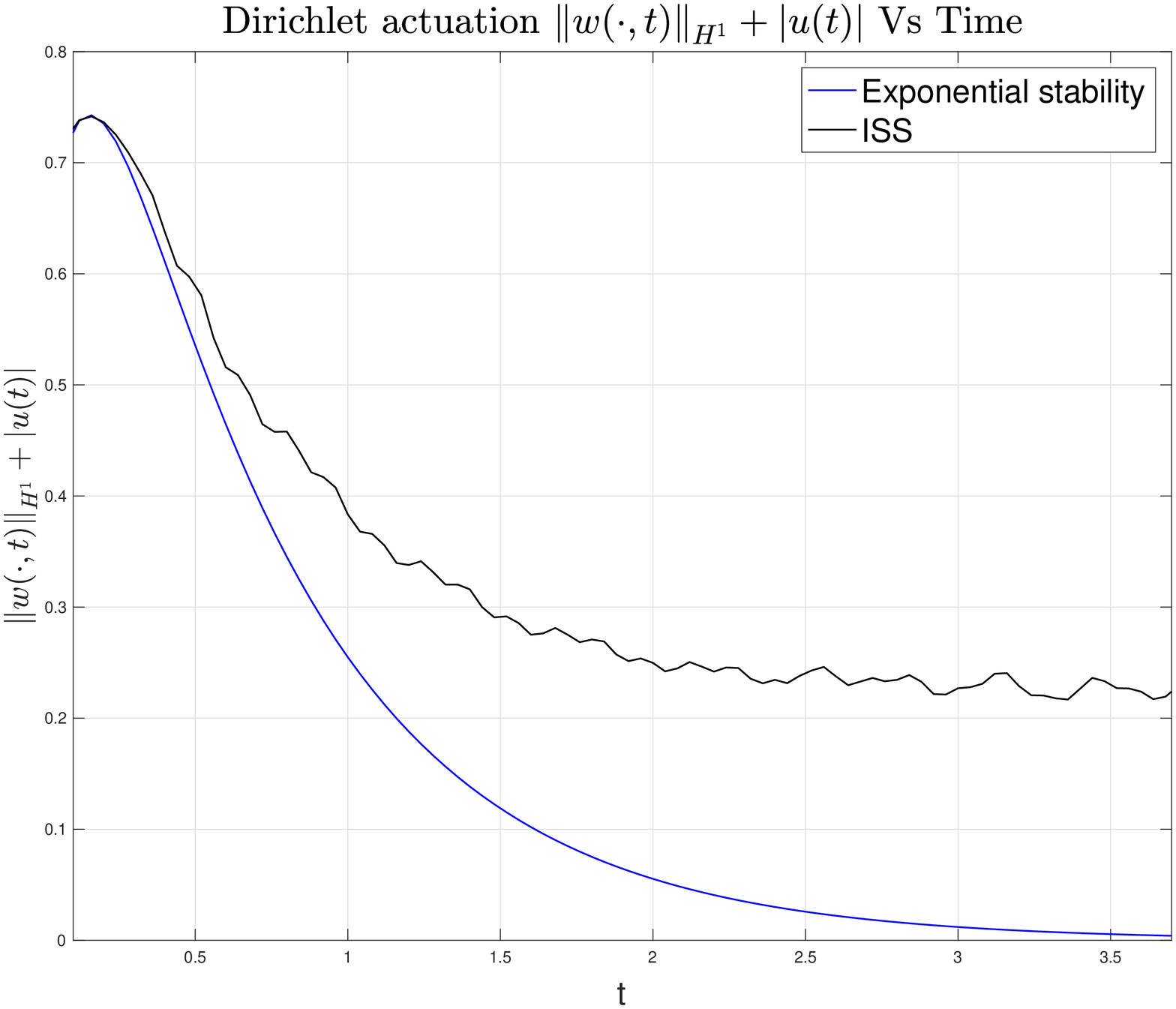}
	\caption{Dirichlet actuation: $\left\|w(\cdot,t)\right\|_{H^1}+\left|u(t)\right|$ vs. $t$}\label{fig:Fig1}
\end{figure}

Next, we carry out two simulations of the closed-loop system for the unperturbed and perturbed cases. In both simulations we have 
$N=4$ and gains given by \eqref{eq:GainsDirStab}. We choose initial conditions
\begin{equation}\label{eq:DirICs}
\begin{array}{lll}
& u(0) = 0,\ \ z(x,0) = w(x,0) = 25(x-x^2)^3, \ x\in [0,1].
\end{array}
\end{equation}
Note that $w(\cdot,0)\in \mathcal{D}\left(\mathcal{A}\right)$, where $\mathcal{D}\left(\mathcal{A}\right)$ is defined in \eqref{eq:D1}. The $H^1$ norm of $w(\cdot,t)$ is approximated by truncating \eqref{lem22} after $60$ coefficients. Then, the ODEs \eqref{eq:WOdesNonDelayedHInf} with $1\leq n \leq 60$ and \eqref{eq:WobsODENonDelayed} are simulated using MATLAB with $v(t) = K_0\hat{w}^{N_0}(t)$ and $\hat{w}^{N_0}(t)$ defined in \eqref{eq:WContDef}. The value of $\zeta_N(t)$ in \eqref{eq:zetaintegral} is approximated using
\begin{equation}\label{eq:zetaNApprox}
\zeta_N(t) \approx \sum_{n=5}^{60} w_n(t)\phi_n^D(x_*).
\end{equation}
In the perturbed case, we consider the disturbances
\begin{equation}\label{eq:SimDist}
d(x,t)=0.25\sin(10x+t),\ \sigma(t)=0.25\cos(30t).
\end{equation}
The simulation results are presented in Figure \ref{fig:Fig1}. From the simulations of exponential stability, we
obtain a decay rate $1.17$, which is slightly larger than the theoretical decay rate $\delta = 1$ found from the LMIs.

\emph{B. Neumann actuation and collocated measurement}\label{ex:KSENeum}

Consider the unperturbed KSE \eqref{eq:LinearizedKSE}, boundary conditions \eqref{eq:BCSNe} and unperturbed measurement \eqref{eq:BoundMeas}. We choose $\delta=1$, which results in $N_0=1$. {To guarantee minimal value of $N$, the gains $K_0$ and $L_0$ were found by solving \eqref{eq:GainsDesignLNeumann} and \eqref{eq:GainsDesignKNeumann} with strong inequality replaced by equality and $\delta=\delta_0=5$. The obtained observer and controller gains are}
\begin{equation}\label{eq:GainsNeuStab}
K_0 = \left[477.83, 32.61, -3315.44 \right], \ L_0=\left[-6.147,8.101\right]^T.
\end{equation}
The LMI of Theorem \ref{Thm:WdynExtensionNeumann} is feasible for  $\Gamma=1$ and minimal $N=6$. Next, we consider the perturbed KSE \eqref{eq:LinearizedKSEHinf}, boundary conditions \eqref{eq:BCSNe} and perturbed measurement \eqref{eq:InDomPointMeasHinfNeu}.
The disturbances again satisfy \eqref{eq:AssumpDisturb}.
For the case of ISS we choose $K_0$ and $L_0$ given by \eqref{eq:GainsNeuStab}.
For $L^2$-gain analysis we consider $\rho_w= 0.1$ and $\rho_u = 0.2$ $\delta=0$. {The gains $K_0$ and $L_0$ were found by solving \eqref{eq:GainsDesignLNeumann} and \eqref{eq:GainsDesignKNeumann} with strong inequality replaced by equality and $\delta=\delta_0=1$. The corresponding gains are}
\begin{equation}\label{eq:GainsHInfNeumEx}
K_0 = \left[291.602, 13.311, -2043.3 \right], \ L_0=\left[-1.967,3.741\right]^T.
\end{equation}
Let $\Gamma = 1$. The LMI \eqref{eq:LMIsHInfNeu1} (with $\delta=0$ and gains \eqref{eq:GainsHInfNeumEx} for $L^2$-gain analysis and with $\delta=1$ and gains \eqref{eq:GainsNeuStab} for ISS) was verified for $N \in \left\{5,7,9,11,13\right\}$. For each choice of $N$, we find the smallest $\gamma$ which guarantees the feasibility of the LMI. The results are presented in Table \ref{Tab:SimISSNeu}. Also in this case, for ISS $\gamma$ decreases as $N$ grows, whereas for $L^2$-gain the resulting $\gamma$ does not grow for larger $N$.

\begin{figure}
	\centering
	\includegraphics[width=50mm,scale=0.11]{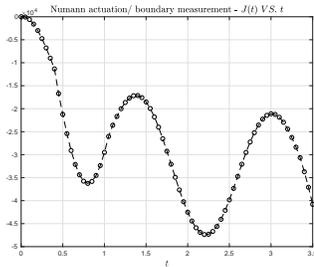}
	\caption{Neumann actuation: $J(t)$ v.s $t$}\label{fig:Fig2}
\end{figure}

\begin{table}
	\begin{center}
		\scalebox{0.70}{
			\begin{tabular}{|l|l|l|l|l|l|l|l|}
				\hline
				N & 5 & 7 & 9 &11 & 13 \\
				\hline
				$\gamma$ (ISS) &3.6 & 1.7 & 1 & 0.6 & 0.5 \\
				\hline
				$\gamma$ ($L^2\text{-gain}$) &31 & 31 & 31 & 31 & 31 \\
				\hline
		\end{tabular}}
	\end{center}
	\caption{\label{Tab:SimISSNeu} Feasibility of LMIs - Neumann. $N$ vs minimal $\gamma$. }
\end{table}
Next, we perform a simulation for the corresponding $L^2$-gain with $\gamma = 31$ and $N=5$. The observer and controller gains are given by \eqref{eq:GainsHInfNeumEx}. The chosen disturbances are given by \eqref{eq:SimDist}. We choose zero initial conditions.
For $t\in [0,3.5]$ we simulate the ODEs \eqref{eq:WOdesNonDelayedHInfNeumann}, $0\leq n \leq 60$ and \eqref{eq:WobsODENonDelayedNeumann} with $v(t)$ defined in \eqref{eq:WContDefNeumann}. The value of $\zeta_N(t)$ in \eqref{eq:zetaintegral} is approximated similarly to \eqref{eq:zetaNApprox}. By truncating Parseval's equality at $n=60$ we approximate the value of
\begin{equation*}
\begin{array}{lll}
&J(t)=\int_0^{t} \left[\rho_w^2\left\|w(\cdot,\tau)\right\|_{L^2}^2+\rho_u^2u^2(\tau)\right.\\
&\hspace{25mm}\left.-\gamma^2\left(\left\|d(\cdot,\tau)\right\|_{L^2}^2 +\sigma^2(\tau)\right) \right]d\tau.
\end{array}
\end{equation*}
Results appear in Figure \ref{fig:Fig2}, confirming the theoretical analysis. We also carry out simulations with $\gamma$ less than $31$ (obtained in LMIs). Simulations show that it is possible to reduce $\gamma$ to approximately $18$, while maintaining $J(t)\leq 0$ for $t\in [0,3.5]$. The latter may indicate the conservatism of the LMIs.

\section{Conclusions}

This paper introduced finite-dimensional observer-based  boundary controllers for linear parabolic PDEs under point measurement. For the 1D linear KSE, modal decomposition using eigenfunctions of a Sturm-Liouville operator and  dynamic extension, with the direct Lyapunov method led to  easily verifiable LMIs for finding the observer dimension. The results were presented for stabilization with guaranteed $L^2$-gain and ISS gain. The presented method allows for challenging finite-dimensional observer-based control of various PDEs, and for design in the case of delayed  inputs and outputs.

\bibliographystyle{IEEEtran}
\bibliography{IEEEabrv,BibliographyKSE}

\begin{thebibliography}{10}
\providecommand{\url}[1]{#1}
\csname url@samestyle\endcsname
\providecommand{\newblock}{\relax}
\providecommand{\bibinfo}[2]{#2}
\providecommand{\BIBentrySTDinterwordspacing}{\spaceskip=0pt\relax}
\providecommand{\BIBentryALTinterwordstretchfactor}{4}
\providecommand{\BIBentryALTinterwordspacing}{\spaceskip=\fontdimen2\font plus
\BIBentryALTinterwordstretchfactor\fontdimen3\font minus
  \fontdimen4\font\relax}
\providecommand{\BIBforeignlanguage}[2]{{%
\expandafter\ifx\csname l@#1\endcsname\relax
\typeout{** WARNING: IEEEtran.bst: No hyphenation pattern has been}%
\typeout{** loaded for the language `#1'. Using the pattern for}%
\typeout{** the default language instead.}%
\else
\language=\csname l@#1\endcsname
\fi
#2}}
\providecommand{\BIBdecl}{\relax}
\BIBdecl

\bibitem{christofides2001}
P.~Christofides, \emph{Nonlinear and Robust Control of PDE Systems: Methods and
  Applications to transport reaction processes}.\hskip 1em plus 0.5em minus
  0.4em\relax Springer, 2001.

\bibitem{kuramoto1975formation}
Y.~Kuramoto and T.~Tsuzuki, ``On the formation of dissipative structures in
  reaction-diffusion systems: Reductive perturbation approach,'' \emph{Progress
  of Theoretical Physics}, vol.~54, no.~3, pp. 687--699, 1975.

\bibitem{sivashinsky1977nonlinear}
G.~Sivashinsky, ``Nonlinear analysis of hydrodynamic instability in laminar
  flames--{I. Derivation} of basic equations,'' \emph{Acta astronautica},
  vol.~4, pp. 1177--1206, 1977.

\bibitem{nicolaenko1986some}
B.~Nicolaenko, ``Some mathematical aspects of flame chaos and flame
  multiplicity,'' \emph{Physica D: Nonlinear Phenomena}, vol.~20, no.~1, pp.
  109--121, 1986.

\bibitem{armaou2000feedback}
A.~Armaou and P.~D. Christofides, ``Feedback control of the
  \uppercase{K}uramoto-\uppercase{S}ivashinsky equation,'' \emph{Physica D:
  Nonlinear Phenomena}, vol. 137, no. 1-2, pp. 49--61, 2000.

\bibitem{christofides2000global}
P.~D. Christofides and A.~Armaou, ``Global stabilization of the
  \uppercase{K}uramoto-\uppercase{S}ivashinsky equation via distributed output
  feedback control,'' \emph{Systems \& Control Letters}, vol.~39, no.~4, pp.
  283--294, 2000.

\bibitem{liu2001stability}
W.-J. Liu and M.~Krsti{\'c}, ``Stability enhancement by boundary control in the
  \uppercase{K}uramoto-\uppercase{S}ivashinsky equation,'' \emph{Nonlinear
  Analysis: Theory, Methods \& Applications}, vol.~43, no.~4, pp. 485--507,
  2001.

\bibitem{cerpa2010null}
E.~Cerpa, ``Null controllability and stabilization of the linear
  \uppercase{K}uramoto-\uppercase{S}ivashinsky equation,'' \emph{Commun. Pure
  Appl. Anal}, vol.~9, no.~1, pp. 91--102, 2010.

\bibitem{guzman2019stabilization}
P.~Guzm{\'a}n, S.~Marx, and E.~Cerpa, ``Stabilization of the linear
  \uppercase{K}uramoto-\uppercase{S}ivashinsky equation with a delayed boundary
  control,'' \emph{IFAC-PapersOnLine}, vol.~52, no.~2, pp. 70--75, 2019.

\bibitem{cerpa2017control}
E.~Cerpa, P.~Guzm{\'a}n, and A.~Mercado, ``On the control of the linear
  {K}uramoto-{S}ivashinsky equation,'' \emph{ESAIM: Control, Optimisation and
  Calculus of Variations}, vol.~23, no.~1, pp. 165--194, 2017.

\bibitem{al2018linearized}
R.~Al~Jamal and K.~Morris, ``Linearized stability of partial differential
  equations with application to stabilization of the
  \uppercase{K}uramoto--\uppercase{S}ivashinsky equation,'' \emph{SIAM Journal
  on Control and Optimization}, vol.~56, no.~1, pp. 120--147, 2018.

\bibitem{Aut12}
E.~Fridman and A.~Blighovsky, ``Robust sampled-data control of a class of
  semilinear parabolic systems,'' \emph{Automatica}, vol.~48, pp. 826--836,
  2012.

\bibitem{NetzerAut14}
N.~Bar~Am and E.~Fridman, ``Network-based {$ H_\infty$} filtering of parabolic
  systems,'' \emph{Automatica}, vol.~50, pp. 3139--3146, 2014.

\bibitem{lunasin2017finite}
E.~Lunasin and E.~S. Titi, ``Finite determining parameters feedback control for
  distributed nonlinear dissipative systems-a computational study,''
  \emph{Evolution Equations \& Control Theory}, vol.~6, no.~4, p. 535, 2017.

\bibitem{kang2019distributed}
W.~Kang and E.~Fridman, ``Distributed stabilization of \uppercase{K}orteweg--de
  \uppercase{V}ries--\uppercase{B}urgers equation in the presence of input
  delay,'' \emph{Automatica}, vol. 100, pp. 260--273, 2019.

\bibitem{selivanov2019delayed}
A.~Selivanov and E.~Fridman, ``Delayed {$ H_\infty$} control of 2\uppercase{D}
  diffusion systems under delayed pointlike measurements,'' \emph{Automatica},
  vol. 109, p. 108541, 2019.

\bibitem{curtain1982finite}
R.~Curtain, ``Finite-dimensional compensator design for parabolic distributed
  systems with point sensors and boundary input,'' \emph{IEEE Transactions on
  Automatic Control}, vol.~27, no.~1, pp. 98--104, 1982.

\bibitem{lasiecka2000control}
I.~Lasiecka and R.~Triggiani, \emph{Control theory for partial differential
  equations: Volume 1, Abstract parabolic systems: Continuous and approximation
  theories}.\hskip 1em plus 0.5em minus 0.4em\relax Cambridge University Press,
  2000, vol.~1.

\bibitem{orlov2004robust}
Y.~Orlov, Y.~Lou, and P.~D. Christofides, ``Robust stabilization of
  infinite-dimensional systems using sliding-mode output feedback control,''
  \emph{International Journal of Control}, vol.~77, no.~12, pp. 1115--1136,
  2004.

\bibitem{katz2020boundary}
R.~Katz, E.~Fridman, and A.~Selivanov, ``Boundary delayed observer-controller
  design for reaction-diffusion systems,'' \emph{IEEE Transactions on Automatic
  Control}, 2021.

\bibitem{balas1988finite}
M.~J. Balas, ``Finite-dimensional controllers for linear distributed parameter
  systems: exponential stability using residual mode filters,'' \emph{Journal
  of Mathematical Analysis and Applications}, vol. 133, no.~2, pp. 283--296,
  1988.

\bibitem{harkort2011finite}
C.~Harkort and J.~Deutscher, ``Finite-dimensional observer-based control of
  linear distributed parameter systems using cascaded output observers,''
  \emph{International journal of control}, vol.~84, no.~1, pp. 107--122, 2011.

\bibitem{selivanov2018boundary}
A.~Selivanov and E.~Fridman, ``Boundary observers for a reaction-diffusion
  system under time-delayed and sampled-data measurements,'' \emph{IEEE
  Transactions on Automatic Control}, vol.~64, no.~4, pp. 3385--3390, 2019.

\bibitem{RamiContructiveFiniteDim}
R.~Katz and E.~Fridman, ``Constructive method for finite-dimensional
  observer-based control of 1-\uppercase{D} parabolic \uppercase{PDE}s,''
  \emph{Automatica}, vol. 122, p. 109285, 2020.

\bibitem{RamiContructiveFiniteDimDelay}
------, ``Delayed finite-dimensional observer-based control of 1-\uppercase{D}
  parabolic \uppercase{PDE}s,'' \emph{Automatica}, vol. 123, p. 109364, 2021.

\bibitem{van2012h}
B.~Van~Keulen, \emph{{$ H_\infty$}-control for distributed parameter systems: A
  state-space approach}.\hskip 1em plus 0.5em minus 0.4em\relax Springer
  Science \& Business Media, 2012.

\bibitem{full}
E.~Fridman and U.~Shaked, ``A descriptor system approach to {$H_{\infty}$}
  control of linear time-delay systems,'' \emph{IEEE Transactions on Automatic
  control}, vol.~47, no.~2, pp. 253--270, 2002.

\bibitem{Aut09b}
E.~Fridman and Y.~Orlov, ``An {LMI} approach to {$H_{\infty}$} boundary control
  of semilinear parabolic and hyperbolic systems,'' \emph{Automatica}, vol.~45,
  no.~9, pp. 2060--2066, 2009.

\bibitem{Karafyllis2016a}
I.~Karafyllis and M.~Krstic, ``{ISS With Respect To Boundary Disturbances for
  1-D Parabolic PDEs},'' \emph{IEEE Transactions on Automatic Control},
  vol.~61, no.~12, pp. 1--23, 2016.

\bibitem{lhachemi2020exponential}
H.~Lhachemi, R.~Shorten, and C.~Prieur, ``Exponential input-to-state
  stabilization of a class of diagonal boundary control systems with delay
  boundary control,'' \emph{Systems \& Control Letters}, vol. 138, p. 104651,
  2020.

\bibitem{jacob2019non}
B.~Jacob, A.~Mironchenko, J.~R. Partington, and F.~Wirth, ``Non-coercive
  {Lyapunov} functions for input-to-state stability of infinite-dimensional
  systems,'' \emph{SIAM Journal on Control and Optimization}, vol.~58, no.~5,
  pp. 2952--2978, 2020.

\bibitem{mironchenko2020input}
A.~Mironchenko and C.~Prieur, ``Input-to-state stability of
  infinite-dimensional systems: recent results and open questions,'' \emph{SIAM
  Review}, vol.~62, no.~3, pp. 529--614, 2020.

\bibitem{curtain}
R.~Curtain and H.~Zwart, \emph{An introduction to infinite-dimensional linear
  systems theory}.\hskip 1em plus 0.5em minus 0.4em\relax Springer, 1995,
  vol.~21.

\bibitem{coron2004global}
J.-M. Coron and E.~Tr{\'e}lat, ``Global steady-state controllability of
  one-dimensional semilinear heat equations,'' \emph{SIAM journal on control
  and optimization}, vol.~43, no.~2, pp. 549--569, 2004.

\bibitem{Rami_CDC20}
R.~Katz and E.~Fridman, ``Finite-dimensional control of the
  {K}uramoto-{S}ivashinsky equation under point measurement and actuation,'' in
  \emph{59th IEEE Conference on Decision and Control}, 2020.

\bibitem{renardy2006introduction}
M.~Renardy and R.~C. Rogers, \emph{An introduction to partial differential
  equations}.\hskip 1em plus 0.5em minus 0.4em\relax Springer Science \&
  Business Media, 2006, vol.~13.

\bibitem{brezis2010functional}
H.~Brezis, \emph{\uppercase{F}unctional analysis, \uppercase{S}obolev spaces
  and partial differential equations}.\hskip 1em plus 0.5em minus 0.4em\relax
  Springer Science \& Business Media, 2010.

\bibitem{zheng2018input}
J.~Zheng and G.~Zhu, ``Input-to-state stability with respect to boundary
  disturbances for a class of semi-linear parabolic equations,''
  \emph{Automatica}, vol.~97, pp. 271--277, 2018.

\bibitem{anders2012higher}
D.~Anders, M.~Dittmann, and K.~Weinberg, ``A higher-order finite element
  approach to the \uppercase{K}uramoto-\uppercase{S}ivashinsky equation,''
  \emph{ZAMM-Journal of Applied Mathematics and Mechanics/Zeitschrift f{\"u}r
  Angewandte Mathematik und Mechanik}, vol.~92, no.~8, pp. 599--607, 2012.

\bibitem{mironchenko2020local}
A.~Mironchenko, C.~Prieur, and F.~Wirth, ``Local stabilization of an unstable
  parabolic equation via saturated controls,'' \emph{IEEE Transactions on
  Automatic Control}, vol.~66, no.~5, pp. 2162--2176, 2020.

\bibitem{GeorgeBook}
M.~Tucsnak and G.~Weiss, \emph{Observation and control for operator
  semigroups}.\hskip 1em plus 0.5em minus 0.4em\relax Springer, 2009.

\bibitem{pazy1983semigroups}
A.~Pazy, \emph{Semigroups of linear operators and applications to partial
  differential equations}.\hskip 1em plus 0.5em minus 0.4em\relax Springer New
  York, 1983, vol.~44.

\bibitem{Fridman14_TDS}
E.~Fridman, \emph{Introduction to time-delay systems: analysis and
  control}.\hskip 1em plus 0.5em minus 0.4em\relax Birkhauser, Systems and
  Control: Foundations and Applications, 2014.

\bibitem{karafyllis2018sampled}
I.~Karafyllis and M.~Krstic, ``Sampled-data boundary feedback control of
  1-\uppercase{D} parabolic \uppercase{PDE}s,'' \emph{Automatica}, vol.~87, pp.
  226--237, 2018.

\end{thebibliography}
\end{document}